\documentclass[12pt]{article} 

\usepackage{amsmath,amssymb,amstext}
\usepackage{color}
\usepackage{algorithm}
\usepackage{algpseudocode}
\usepackage{enumerate}
\usepackage{geometry} 
 \numberwithin{equation}{section}

\title{Convergence of cyclic coordinatewise $\ell_1$ minimization}
\author{Kshitij Khare and Bala Rajaratnam\\University of Florida and Stanford University}
\date{}

\begin{document}
\maketitle

\begin{abstract}
We consider the general problem of minimizing an objective function which is the sum of a convex 
function (not strictly convex) and absolute values of a subset of variables (or equivalently the 
$\ell_1$-norm of the variables). This problem appears extensively in modern statistical 
applications associated with high-dimensional data or ``big data", and corresponds to optimizing 
$\ell_1$-regularized likelihoods in the context of model selection. In such applications, cyclic 
coordinatewise minimization (CCM), where 
the objective function is sequentially minimized with respect to each individual coordinate, is often 
employed as it offers a computationally cheap and effective optimization method. Consequently, it is crucial to obtain theoretical guarantees of convergence for the sequence of iterates produced by the cyclic coordinatewise minimization in this setting. Moreover, as the objective corresponds to flat $\ell_1$-regularized likelihoods of many variables, it is important to obtain convergence of the iterates themselves, and not just the function values. Previous results in the literature only establish either, (i) that every limit point of the sequence of iterates is a stationary point of the objective function, or (ii) establish convergence under special assumptions, or (iii) establish convergence for a different minimization approach (which uses quadratic approximation based gradient descent followed by an inexact line search), (iv) establish convergence of only the function values of the sequence of iterates produced by random coordinatewise minimization (a variant of CCM). In this paper, a rigorous general proof of convergence for the cyclic coordinatewise minimization algorithm is provided. We demonstrate the usefulness of our general results in contemporary applications by employing them to prove convergence of two algorithms commonly used in high-dimensional covariance estimation and logistic regression. 
\end{abstract}

\section{Introduction}

\noindent
Let $g: \mathbb{R}^m \rightarrow \mathbb{R} \cup \{\infty\}$ be a twice differentiable strictly convex 
function, whose effective domain has a non-empty interior $C_g$. Suppose also that $g$ has a 
positive curvature everywhere on $C_g$, and that $g({\bf t})$ converges to infinity as ${\bf t}$ 
approaches the boundary of $C_g$. Let $S$ be a given subset of $\{1,2, \cdots, n\}$, $E$ be an 
$m\times n$ matrix having no zero column, and $\lambda > 0$ be fixed. Let 
$$
\mathcal{X} = \{{\bf x} \in \mathbb{R}^n: \; x_i \geq 0 \mbox{ for every } i \in S^c\}. 
$$

\noindent
Define the functions $f_1$ and $f_2$, where $f_i : \mathbb{R}^n \rightarrow \mathbb{R}$ for $i = 1,2$ as follows: 
$$
f_1 ({\bf x}) = g(E{\bf x}) + \lambda \sum_{i \in S} |x_i|, 
$$

\noindent
and, 
$$
f_2 ({\bf x}) = {\bf x}^T E^T E {\bf x} - \sum_{i \in S^c} \log x_i + \lambda \sum_{i \in S} |x_i|. 
$$

\noindent
Consider the following two minimization problems: 
\begin{equation} \label{l1minpblm1}
\mbox{Minimize } f_1 ({\bf x}) \mbox{ subject to } {\bf x} \in \mathcal{X}. 
\end{equation}

\noindent
\begin{equation} \label{l1minpblm2}
\mbox{Minimize } f_2 ({\bf x}) \mbox{ subject to } {\bf x} \in \mathcal{X}. 
\end{equation}

\noindent
The minimization problems in (\ref{l1minpblm1}) and (\ref{l1minpblm2}) appear extensively in 
contemporary applications, and are particularly relevant in statistics and machine learning (see 
for example \cite{Friedman:Hastie:2007, FHT:2010, Fu:1998, KOR:2014:CONCORD, 
Meinshausen:Buhlmann:2006, Shevade:Keerthi:2003, Tibshirani:1996, Tibshirani:Saunders:2003}), and signal processing (see for example \cite{Bradley:Fayyad:1999, Chen:Donoho:1999, Donoho:Johnstone:1994, Donoho:Johnstone:1995, Sardy:Bruce:2001, Sardy:Tseng:2004}). In statistical applications, the function $g(E {\bf x})$ is typically a log-likelihood 
or pseudo log-likelihood corresponding to a statistical model. Traditional statistical methods focus on 
minimizing the function $g(E {\bf x})$ without the addition of a minimizer like $\|{\bf x}\|_1$. 
However, in the modern context of high-dimensional data or ``big data", it is often desirable to obtain sparse solutions (solutions where many entries are exactly equal to zero), resulting in the inclusion 
of the term $\lambda \sum_{i \in S} |x_i|$ in the objective function. The most challenging features of 
the minimization problems (\ref{l1minpblm1}) and (\ref{l1minpblm2}) are the following: 
\begin{enumerate}
\item In many applications, $m < n$. Hence the functions $g(E {\bf x})$, $f_1 ({\bf x})$ and $f_2 ({\bf x})$ are not 
strictly convex, and in general do not have a unique global minimum. 
\item The minimization occurs on a high-dimensional space, with hundreds or thousands (if not more) 
of variables. 
\item The minimization problem is non-smooth due to the presence of the ``$\ell_1$ penalty" 
term\\ 
$\lambda \sum_{i \in S} |x_i|$. 
\end{enumerate}

\noindent
Hence, any method proposed for finding a solution to the above problem should be computationally 
scalable and have theoretical convergence guarantees. In many statistical applications involving 
high dimensional regression and high dimensional covariance estimation (see for example 
\cite{Fu:1998, FHHT:2007, KOR:2014:CONCORD}), coordinatewise minimization can be 
performed in closed form. Hence, for such problems, a cyclic coordinatewise minimization (CCM) 
algorithm (where each iteration consists of minimizing the objective function sequentially 
over all the coordinates) is often used, as it offers a computationally cheap and effective method for 
minimizing the respective objective functions. In situations where coordinatewise minimization 
cannot be achieved in closed form, it often involves minimizing a one-dimensional convex function, 
and can be numerically achieved to a high degree of accuracy in a few steps. Hence, coordinatewise 
minimization has also been used in such situations (see for example \cite{Shevade:Keerthi:2003}). 
Hence, understanding the convergence properties of the cyclic coordinatewise minimization 
algorithm for (\ref{l1minpblm1}) and (\ref{l1minpblm2}) is a crucial and relevant task. However, a 
rigorous proof of convergence of the cyclic coordinatewise minimization algorithm for 
minimization problems in (\ref{l1minpblm1}) and (\ref{l1minpblm2}) is not available in the 
literature. We now provide a brief overview of existing optimization methods and convergence 
results related to these problems.  

In Tseng \cite{Tseng:2001}, it is proved that under appropriate conditions on $g$, every limit point of 
the sequence of iterates produced by the cyclic coordinatewise minimization algorithm is a stationary 
point of the corresponding objective function. However, this does not necessarily mean that the 
sequence of iterates converges. Tseng and Yun \cite{Tseng:Yun:2009} propose a 
block-coordinatewise gradient descent (CGD) approach, which can be thought of as a hybrid of 
gradient-projection and coordinate descent. In particular, they consider minimizing an objective 
function of the form 
\begin{equation} \label{tsengyunmp}
f({\bf x}) + c \sum_{i=1}^n P_j (x_j), 
\end{equation}

\noindent
where $P_j$ is a proper, convex, lower semicontinuous function for $1 \leq j \leq n$, and $f$ 
is a continuously differentiable function on an open subset of $\mathbb{R}^n$ containing the 
effective domain of $P_j$ for every $1 \leq j \leq n$. At each iteration, a quadratic approximation 
of the function $f$ is considered, a descent direction is then generated by applying block 
coordinate descent, followed by an inexact line search along this direction (by using an 
Armijo-type rule to ensure sufficient descent). The authors in \cite{Tseng:Yun:2009} also provide a proof that the sequence 
of iterates produced by their algorithm converges under suitable assumptions. Note that if the 
function $g$ in (\ref{l1minpblm1}) is not quadratic, then clearly the CCM and CGD approaches are 
distinctly different. Considering (\ref{l1minpblm2}), we note that it can be expressed in the 
framework of (\ref{tsengyunmp}) in two ways. We can choose $f({\bf x}) = {\bf x}^T E^T E {\bf x} 
- \sum_{i \in S^c} \log x_i$ and $P_j (x_j) = |x_j|$ for $j \in S$, in which case the CCM 
and CGD approaches are again different as $f$ is not quadratic. Alternatively, if 
we choose $f({\bf x}) = {\bf x}^T E^T E {\bf x}$, $P_j (x_j) = |x_j|$ for $j \in S$, and $P_j (x_j) = 
- \log x_j$ for $j \in S^c$, then the function $\sum_{j=1}^n P_j (x_j)$ is not polyhedral. Hence, 
the assumptions in \cite[Lemma 7]{Tseng:Yun:2009} do not apply, and it is not clear if the 
convergence results in \cite[Theorem 2]{Tseng:Yun:2009} and \cite[Theorem 3]{Tseng:Yun:2009} 
apply. 

Saha and Tewari \cite{Saha:Tewari:2010} provide finite time convergence results for a 
variety of cyclic coordinatewise descent methods for objective functions of the form $f({\bf x}) + 
\lambda \sum_{j=1}^n |x_j|$. However, their convergence results rely on the assumption that the 
function $f$ is isotone, i.e., essentially $f$ is twice-differentiable and the Hessian matrix of $f$ at 
any ${\bf x}$ in its effective domain has non-positive off-diagonal entries. Such an assumption 
does not hold in general for many contemporary applications and those which we consider in 
Section \ref{sec:applications}. 

Luo and Tseng \cite{Luo:Tseng:1989} consider the following minimization problem. 
\begin{equation} \label{luotsengmp}
\mbox{Minimize } g(E{\bf x}) + {\bf b}^T {\bf x} \mbox{ subject to } l_i \leq x_i \leq u_i \; \forall \; 
1 \leq i \leq n, 
\end{equation}

\noindent
where ${\bf b}, {\bf l}, {\bf u}$ are fixed $n$-dimensional vectors. The entries of ${\bf l}$ and ${\bf 
u}$ are allowed to be $-\infty$ and $\infty$ respectively. They provide a very detailed and intricate 
proof of convergence of the sequence of iterates produced by the cyclic coordinatewise descent 
algorithm for (\ref{luotsengmp}) (see also \cite{Luo:Tseng:1:1989, Luo:Tseng:1992}). Note once more that 
the minimization problem in (\ref{luotsengmp}) is substantially different than the minimization 
problem in (\ref{l1minpblm1}) and (\ref{l1minpblm2}). 

In recent useful work, Richtarik and Takac \cite{Richtarik:Takac:2011} provide a random coordinatewise descent algorithm for solving (\ref{tsengyunmp}), where instead of cycling over all the coordinate blocks, a random coordinate is chosen and minimized over at each iteration. Intuitively speaking, a randomized choice of coordinates may avoid a possible worst case ordering of the coordinates in the cyclic setting, is also more suitable for 
situations when all the data in not available all the time, and more amenable for a convergence analysis. The authors in 
\cite{Richtarik:Takac:2011} establish  important convergence (and provides rates) for the function values of the sequence of iterates produced by the random coordinatewise descent algorithm. Establishing convergence of  the sequence of \emph{iterates} for random coordinatewise descent however remains a challenge. We note that one of the compelling reasons that has motivated the use of random coordinatewise descent algorithm versus the (non-random) cyclic coordinatewise minimization is that the former allows for easier convergence analysis, though many methods that have been proposed in the machine learning and statistics literature actually use the (non-random) cyclic coordinatewise minimization. In this paper we address the crucial and challenging problem of establishing convergence of the sequence of iterates in the (non-random) cyclic coordinatewise minimization setting. We also note that in modern high-dimensional problems in Statistics/Machine Learning, the objective functions are often very flat, and it is quite likely that although the function values converge, the sequence of iterates do not. 

Several methods other than cyclic coordinatewise minimization have also been proposed in the literature to 
solve the minimization problems in (\ref{l1minpblm1}) and (\ref{l1minpblm2}) (or a more general 
version of this problem, where the term $\lambda \sum_{i \in S} |x_i|$ is replaced by a (block) 
separable non-smooth function). One class of methods is based on proximal gradient descent with 
an Armijo-type stepsize (see for example \cite{Fukushima:Mine, Kiwiel}). Another class of methods is 
based on trust-regions (see for example \cite{Auslender:1981, Burke:1985, Fletcher:1982}). See 
Tseng and Yun \cite{Tseng:Yun:2009} for a detailed list of related references. We note that none of 
these methods correspond to the classical coordinatewise minimization approach that has been proposed and 
extensively used in the statistical applications outlined above. 

In this paper, we provide a rigorous proof of convergence of the sequence of iterates produced by the 
cyclic coordinatewise minimization algorithm for the minimization problems in (\ref{l1minpblm1}) and 
(\ref{l1minpblm2}). We shall build on the work of Luo and Tseng \cite{Luo:Tseng:1989}, and extend it 
when incorporating non-differentiable terms of the form $\lambda \sum_{i \in S} |x_i|$ in $f_1({\bf x})$ and 
$f_2 ({\bf x})$. This generalization makes the convergence analysis of the cyclic coordinatewse minimization 
algorithm for (\ref{l1minpblm1}) and (\ref{l1minpblm2}) more complex as compared to the convergence 
analysis of the cyclic coordinatewise minimization approach for (\ref{luotsengmp}). We shall see that the 
non-smooth term leads to many challenging and non-trivial questions. 

This paper is organized as follows. In Section \ref{section:summary:results}, we provide a 
summary of the assumptions, algorithms and the main convergence results in the paper. A detailed 
proof of convergence for the cyclic coordinatewise minimization algorithm for the minimization 
problems in (\ref{l1minpblm1}) and (\ref{l1minpblm2}) is then provided in Section \ref{sec:convergence} 
and Section \ref{sec:convergence:withlog} respectively. The results in Section \ref{sec:convergence} and 
Section \ref{sec:convergence:withlog} are then used in Section \ref{sec:applications} to establish 
convergence of two algorithms arising in high-dimensional covariance estimation and high dimensional 
logistic regression.

\section{Summary of main results} \label{section:summary:results}

\noindent
In this section, we undertake the following: (a) provide the assumptions that are made for the minimization problems in 
(\ref{l1minpblm1}) and (\ref{l1minpblm2}), (b) formally define the cyclic coordinatewise 
minimization algorithms corresponding to these problems, and (c) state the main convergence 
results that are established later in this paper. Recall that $f_1 ({\bf x}) = g(E{\bf x}) + 
\lambda \sum_{i \in S} |x_i|$. We start by providing the assumptions that will be made for the 
minimization problem in 
(\ref{l1minpblm1}). 
\begin{itemize}
\item (A1) The effective domain of $g$ has a non-empty interior $C_g$. 
\item (A2) $g$ is strictly convex and twice continuously differentiable on $C_g$. 
\item (A3) Either $g(t) \rightarrow \infty$ as $t$ approaches the boundary of $C_g$, or, 
$|S| = n$ and $g$ is non-negative with $C_g = \mathbb{R}^m$. 
\item (A4) $g$ has a positive curvature everywhere on $C_g$. 
\item (A5) The set of optimal solutions of the minimization problem in (\ref{l1minpblm1}), 
denoted by $\mathcal{X}^*$, is non-empty. 
\end{itemize}

\noindent
Consider the following practical implementation of the coordinatewise descent (CCM) algorithm to solve the minimization 
problem in (\ref{l1minpblm1}). 
\begin{figure}[H]
  \centering
  \begin{minipage}[t]{0.9\textwidth}
    \alglanguage{pseudocode}
    \begin{algorithm}[H]
      \caption{Cyclic coordinatewise descent algorithm for $f_1$}
      \begin{algorithmic}
        \State 1. Set $r = 0$. Start with initial value ${\bf x}^0 \in \mathcal{X}$ such that $f_1 ({\bf x}^0)$ 
is finite, and a prespecified tolerance $\epsilon$. 
        \State 2. Set ${\bf x}^{r,0} = {\bf x}^r$. 
        \State 3. For $i = 1,2, \cdots, n$, set 
\begin{equation} \label{eq3}
{\bf x}^{r,i} = arg min_{{\bf x} \in \mathcal{X}: \; x_j = x^{r,i-1}_j \forall j \neq i} f_1 ({\bf x}). 
\end{equation}
        \State 4. Set ${\bf x}^{r+1} = {\bf x}^{r,n}$. If $\|{\bf x}^{r+1} - {\bf x}^r\| > \epsilon$, 
set r = r + 1, return to Step $2$. Otherwise, stop.  
      \end{algorithmic}
    \end{algorithm}
  \end{minipage}  
\end{figure}

\noindent
We first claim that (\ref{eq3}) is well-defined by using contradiction.  Note that for any $\xi \in 
\mathbb{R}$, the set $H_\xi := \{E{\bf x}: \; {\bf x} \in \mathcal{X}, f_1 ({\bf x}) \leq \xi\}$ is 
contained in the set $\{E{\bf x}: \; {\bf x} \in \mathcal{X}, g(E {\bf x}) \leq \xi\}$. It follows by 
\cite[Lemma A.1]{Luo:Tseng:1989} that if $g(t) \rightarrow \infty$ as $t$ approaches the boundary of 
$C_g$, then $\{E{\bf x}: \; {\bf x} \in \mathcal{X}, g(E {\bf x}) \leq \xi\}$ is bounded. Alternatively, if 
$|S| = n$ and $g$ is non-negative, then $H_\xi$ is contained in the set 
$\{E{\bf x}: \; \sum_{i=1}^n |x_i| \leq \xi/\lambda\}$. In either case, we get that 
\begin{equation} \label{eq6.1}
H_\xi \mbox{ is bounded for every } \xi \in \mathbb{R}. 
\end{equation}

\noindent
Suppose now that the minimum in (\ref{eq3}) is not attained for some $r$ and $i$. Let ${\bf e}^i$ 
denote the $i^{th}$ unit vector in $\mathbb{R}^n$. There are then two possibilities: 
\begin{enumerate}[(a)]
\item $i \in S$ and $f_1 ({\bf x}^{r,i-1} - h {\bf e}^i)$ is non-increasing as $h \rightarrow \infty$. 
Hence, ${\bf x}^{r,i-1} - h {\bf e}^i \in H_{f_1 ({\bf x}^{r,i-1})}$ for large enough $h$. The 
boundedness of $H_{f_1 ({\bf x}^{r,i-1})}$ implies that $E {\bf e}^i = 0$, which contradicts the 
assumption that no column of $E$ is zero. 
\item $f_1 ({\bf x}^{r,i-1} + h {\bf e}^i)$ is non-increasing as $h \rightarrow \infty$. This case leads 
to the same contradiction as in (a). 
\end{enumerate}

\noindent
The following theorem now formally establishes convergence of the sequence of iterates produced by 
Algorithm 1, and is the first main result in this paper. 
\newtheorem{thm}{Theorem}[section]
\begin{thm} \label{thm1}
The sequence of iterates $\{{\bf x}^r\}_{r \geq 0}$ generated by the cyclic coordinatewise descent 
algorithm for $f_1$ converges to a value ${\bf x}^* \in \mathcal{X}$ such that $f_1 ({\bf x}^*) \leq 
f_1 ({\bf x})$ for every ${\bf x} \in \mathcal{X}$. 
\end{thm}

\noindent
A detailed proof of Theorem \ref{thm1} will be provided in Section \ref{sec:convergence}. We 
now briefly outline the major steps in the proof. We first show that the difference between the successive 
iterates produced by the CCM algorithm for $f_1$ goes to zero. Further arguments establish that this sequence of 
differences between the successive iterates is actually square-summable. Note that square-summability of the 
sequence of differences is not sufficient to establish that $\{{\bf x}^r\}_{r \geq 0}$ is a Cauchy sequence. We then 
proceed to show that the distance between the sequence of iterates and the boundary of $\mathcal{X}^*$ (the set of optimal 
solutions) goes to zero. Again, this itself is also not sufficient to establish that $\{{\bf x}^r\}_{r \geq 0}$ is a Cauchy sequence 
(see the discussion just before Lemma \ref{lem11}). However, using the three facts above, along with some matrix-theoretic 
results and combinatorial arguments, we prove that the sequence of iterates produced by the cyclic coordinatewise 
descent algorithm is a Cauchy sequence with limit ${\bf x}^* \in \mathcal{X}$. This is done as follows. First, we show 
that eventually some coordinates of ${\bf x}^r$ are exactly equal to zero, while the remaining coordinates are bounded 
from zero as $r \rightarrow \infty$. Second, we show that (see Lemma \ref{lem17}) the coordinates of ${\bf x}^r$ that stay 
away from zero are influenced by those coordinates which eventually become zero. Moreover, this influence is a function of  
the distance between these ultimate zero coordinates and zero, and therefore dies away as $r \rightarrow \infty$. This is 
then used, along with a series of combinatorial arguments, to establish that (see Lemma \ref{lem19}) for an arbitrary 
$\epsilon > 0$, there exists an ${\bf x}^* \in \mathcal{X}^*$ such that $\|{\bf x}^r - {\bf x}^*\| < \epsilon$ for large 
enough $r$. This immediately implies that $\{{\bf x}^r\}_{r \geq 0}$ is a Cauchy sequence. Since we have already 
established that the distance between the sequence of iterates and the boundary of $\mathcal{X}^*$  goes to zero, 
convergence to an optimal solution follows. 

Now, we consider the problem of minimizing the function $f_2$ defined in (\ref{l1minpblm2}). 
Recall that 
$$
{\bf x}^T E^T E {\bf x} - \sum_{i \in S^c} \log x_i + \lambda \sum_{i \in S} |x_i|. 
$$

\noindent
For the function $f_2$, the only assumption that is made is a stronger version of assumption (A5), 
this assumption essentially states that the level sets of $f_2$ are bounded, and is standard in many 
contemporary applications.  
\begin{itemize}
\item (A5)* Let $\xi \in \mathbb{R}$ be arbitrarily fixed. If ${\bf x} \in \mathcal{X}$ satisifes 
$f_2 ({\bf x}) \leq \xi$, then there exists $\xi^* \in \mathbb{R}_+$ (independent of ${\bf x}$) 
such that $1/\xi^* \leq x_i \leq \xi^*$ for every $i \in S^c$ and $|x_i| \leq \xi^*$ for every 
$i \in S$. 
\end{itemize}

\noindent
We shall also show that this level set assumption will se satisfied in the application considered in 
Section \ref{sec:applications}. Again, we consider the following coordinatewise descent 
algorithm to solve the minimization problem in (\ref{l1minpblm2}). Note that the steps of 
the following algorithm are identical to that of Algorithm 1, excpet that $f_1$ is replaced 
by $f_2$. However, we have provided separate statements of the two algorithms for expositional 
convenience, in particular, for differentiating between the sequence of iterates produced by 
applying the CCM algorithm for $f_1$ and $f_2$. 
\begin{figure}
  \centering
  \begin{minipage}[t]{0.9\textwidth}
    \alglanguage{pseudocode}
    \begin{algorithm}[H]
      \caption{Cyclic coordinatewise descent algorithm for $f_2$}
      \begin{algorithmic}
        \State 1. Set $r = 0$. Start with initial value ${\bf z}^0 \in \mathcal{X}$ such that $f_2 ({\bf z}^0)$ 
is finite, and a prespecified tolerance $\epsilon$. 
        \State 2. Set ${\bf z}^{r,0} = {\bf z}^r$. 
        \State 3. For $i = 1,2, \cdots, n$, set 
\begin{equation} \label{eqwlog3}
{\bf z}^{r,i} = arg min_{{\bf x} \in \mathcal{X}, x_j = z^{r,i-1}_j \forall j \neq i} f_2 ({\bf x}). 
\end{equation}
        \State 4. Set ${\bf z}^{r+1} = {\bf z}^{r,n}$. If $\|{\bf z}^{r+1} - {\bf z}^r\| > \epsilon$, 
set r = r + 1, go to Step $2$. 
      \end{algorithmic}
    \end{algorithm}
  \end{minipage}  
\end{figure}

\noindent
It will be shown in Section \ref{sec:convergence:withlog} (see Lemma \ref{lemwlog1}) that the 
minimization in (\ref{eqwlog3}) is well-defined, and the unique minimizer can be obtained in 
closed form. The following theorem establishes convergence of the sequence of iterates 
produced by the cyclic coordinatewise descent algorithm for minimizing $f_2$ and is the second 
main result in this paper. 
\begin{thm} \label{thmwlog1}
The sequence of iterates $\{{\bf z}^r\}_{r \geq 0}$ generated by the cyclic coordinatewise descent 
algorithm for minimizing $f_2$ converges to a ${\bf z}^* \in \mathcal{X}$ such that $f_2 ({\bf z}^*) \leq 
f_2 ({\bf x})$ for every ${\bf x} \in \mathcal{X}$. 
\end{thm}

\noindent
A detailed proof of Theorem \ref{thmwlog1} will be provided in 
Section \ref{sec:convergence:withlog}. There are two differences between the functions $f_1$ and $f_2$. 
Let $q({\bf y}) = {\bf y}^T {\bf y}$ for every ${\bf y} \in \mathbb{R}^m$. The term 
$g(E{\bf x})$ in $f_1$ is replaced by the special choice $q(E {\bf x}) = {\bf x}^T E^T E {\bf x}$ 
in $f_2$. The presence of the logarithmic terms in $f_2$ however introduces a new feature 
as compared to $f_1$. Hence, although the basic method of proving convergence remains the 
same for $f_2$, the presence of the logarithmic terms in $f_2$ create new challenges which 
will be tackled in the convergence analysis in Section \ref{sec:convergence:withlog}.

\section{Convergence analysis for cyclic coordinatewise minimization applied to $f_1$} 
\label{sec:convergence}

\noindent
In this section, we provide a detailed proof of Theorem \ref{thm1}. We start with the following lemma 
about $\mathcal{X}^*$, the set of optimal solutions of the minimization problem (\ref{l1minpblm1}). 
The proof of this lemma follows immediately from arguments in \cite[Page 5]{Luo:Tseng:1989} and 
is therefore omitted. 
\newtheorem{lemma}{Lemma}[section]
\begin{lemma} \label{lem1}
$\mathcal{X}^*$ is a convex set. Also, there exists ${\bf t}^* \in \mathbb{R}^m$ such that 
$$
E {\bf x}^* = {\bf t}^*, \; \forall \; {\bf x}^* \in \mathcal{X}^*. 
$$
\end{lemma}

\noindent
It follows from assumptions (A1) and (A4) that $\nabla^2 g$ is positive definite in some open ball 
$U^*$ containing ${\bf t}^*$. Hence, there exists $\sigma > 0$ such that 
\begin{equation} \label{eq1}
\nabla^2 g ({\bf t}) - \sigma I_n \mbox{ is positive definite } \forall {\bf t} \in U^*. 
\end{equation}

\noindent
Let ${\bf d} ({\bf x}) = \nabla \{g(E {\bf x})\} = E^T \nabla g(E {\bf x})$, where $ \nabla g(E {\bf x})$ 
denotes the gradient function of $g$ evaluated at $E {\bf x}$. We denote the $i^{th}$ 
entry of ${\bf d} ({\bf x})$ by $d_i ({\bf x})$. Let ${\bf d}^* := E^T \nabla g(t^*)$. It follows by 
Lemma \ref{lem1} that 
\begin{equation} \label{eq2}
{\bf d} ({\bf x}^*) = {\bf d}^* \; \forall x^* \in \mathcal{X}^*. 
\end{equation}

\noindent
Note that the sub differential versions of the KKT conditions for the convex minimization problem in (\ref{l1minpblm1}) imply 
that ${\bf x} \in \mathcal{X}^*$ if and only if 
\begin{eqnarray}
& & x_i = \max(0, x_i - d_i ({\bf x}))\; \mbox{ for } i \in S^c, \label{kkt1}\\
& & d_i ({\bf x}) + \lambda sign(x_i) = 0 \; \mbox{ if } x_i \neq 0, i \in S, \label{kkt2}\\
& & |d_i ({\bf x})| \leq \lambda \; \mbox{ if } x_i = 0, i \in S. \label{kkt3} 
\end{eqnarray}

\noindent
We provide an alternative characterization of the elements of $\mathcal{X}^*$, which will be useful 
in our analysis. 
\begin{lemma} \label{lem2}
${\bf x} \in \mathcal{X}^*$ if and only if 
\begin{eqnarray}
& & x_i = \max(0, x_i - d_i ({\bf x}))\; \mbox{ for } i \in S^c, \label{kkt4}\\
& & x_i = sign(x_i - d_i ({\bf x})) \max(|x_i - d_i ({\bf x})| - \lambda, 0) \mbox{ for } i \in S. 
\label{kkt5} 
\end{eqnarray}
\end{lemma}

\noindent
The proof of this lemma is provided in the appendix. Recall that $\{{\bf x}^r\}_{r \geq 0}$ is the sequence of iterates 
generated by the coordinatewise descent algorithm for minimizing $f_1$, and ${\bf x}^{r,i}$ is the appropriate 
coordinatewise minimizer defined in (\ref{eq3}). It follows from the arguments in the proof of Lemma \ref{lem2} that 
for $i \in S$, 
\begin{equation} \label{eq4}
x^{r,i}_i  = sign(x^{r,i}_i - d_i ({\bf x}^{r,i})) \max(|x^{r,i}_i - d_i ({\bf x}^{r,i})| - \lambda, 0), 
\end{equation}

\noindent
and for $i \in S^c$ 
\begin{equation} \label{eq5}
x^{r,i}_i = \max(0, x^{r,i}_i - d_i ({\bf x}^{r,i})). 
\end{equation}

\noindent
Next, we state a lemma from \cite{Hoffman:1952} which was used in \cite{Luo:Tseng:1989}, and will 
also play an important role in our analysis. Let $\|{\bf x}\| := \sqrt{{\bf x}^T {\bf x}}$ denote the 
Euclidean norm, and ${\bf x}^+ := (\max(0, x_i))_{i=1}^n$ for any vector ${\bf x}$. Also, ${\bf x} 
\leq {\bf y}$ implies that $x_i \leq y_i$ for every $1 \leq i \leq n$. 
\begin{lemma}[\cite{Hoffman:1952}] \label{lem3}
Let $B_1$ and $B_2$ be any $k_1 \times n$ and $k_2 \times n$ matrices respectively. Then, there 
exists a constant $\theta > 0$ depending only on $B_1$ and $B_2$ such that, for any ${\bar {\bf x}} 
\in \mathcal{X}$ and any $k_1$-vector ${\bf d}_1$ and $k_2$-vector ${\bf d}_2$ such that the linear 
system $B_1 {\bf y}  = {\bf d}_1, B_2 {\bf y} \leq {\bf d}_2, \; {\bf y} \in \mathcal{X}$ is consistent, there is a point 
${\bar {\bf y}}$ satisfying $B_1 {\bar {\bf y}} = {\bf d}_1, B_2 {\bar {\bf y}} = {\bf d}_2, \; {\bar {\bf y}} 
\in \mathcal{X}$, with 
$$
\|{\bar {\bf x}} - {\bar {\bf y}}\| \leq \theta (\|B_1 {\bar {\bf x}} - {\bf d}_1\| + \|(B_2 {\bar {\bf x}} - 
{\bf d}_2)^+\|). 
$$
\end{lemma}

\noindent
Now let 
$$
{\bf t}^{r,i} = E {\bf x}^{r,i} 
$$

\noindent
for all $r$ and all $0 \leq i \leq n$. By (\ref{eq3}), it follows that 
\begin{equation} \label{eq6} 
f_1 ({\bf x}^{r,i}) \leq f_1 ({\bf x}^{r,i-1}) 
\end{equation}

\noindent 
for every $r$ and $1 \leq i \leq n$. Hence, ${\bf t}^{r,i} \in H_{f_1 ({\bf x}^0)}$ for every $r$ and 
$0 \leq i \leq n$. It follows by (\ref{eq6.1}) that 
\begin{equation} \label{eq6.2}
\{{\bf t}^{r,i}\}_{r \geq 0, 0 \leq i \leq n} \mbox{ is bounded}. 
\end{equation}

\noindent
Also, since $g$ is twice continuously differentiable, it follows that $\{g({\bf t}^{r,i})\}_{r \geq 0}$ is uniformly bounded 
above for all $0 \leq i \leq n$. If $g(t) \rightarrow \infty$ as $t$ approaches the boundary of $C_g$, 
it follows that every limit point of $\{{\bf t}^{r,i}\}_{r \geq 0}$ lies in $C_g$ for all $0 \leq i \leq n$. 
If $C_g = \mathbb{R}^m$, then it follows by (\ref{eq6.2}) that again every limit point of 
$\{{\bf t}^{r,i}\}_{r \geq 0}$ lies in $C_g$ for all $0 \leq i \leq n$. By (\ref{eq6}), the sequence $\{f_1 
({\bf x}^{r,i})\}_{r \geq 0}$ decreases to the same quantity, say $f^\infty$ for every $0 \leq i \leq n$. 
If $f^\infty = -\infty$, then assumption (A5) (which says that the set of optimal solutions to (\ref{l1minpblm1}) is 
non-empty) will be violated. Hence $f^\infty > -\infty$. We now 
prove that the difference between the successive iterates of the cyclic coordinatewise descent algorithm for 
$f_1$ converges to zero. 
\begin{lemma} \label{lem4}
$$
\| {\bf x}^{r+1} - {\bf x}^r\| \rightarrow 0 \mbox{ as } r \rightarrow \infty. 
$$
\end{lemma}

\noindent
{\it Proof} We proceed by contradiction. Suppose the result does not hold. Then there exists 
$\epsilon > 0$, $i \in \{1,2, \cdots, n\}$ and a subsequence $\mathcal{R}$ of $\mathbb{N}$ such 
that $|x_i^{r+1} - x_i^r| > \epsilon$ for every $r \in \mathcal{R}$. It follows by the definition of 
${\bf t}^{r,i}$ that 
\begin{equation} \label{eq7.1}
\|{\bf t}^{r,i} - {\bf t}^{r,i-1}\| = \|E ({\bf x}^{r,i} - {\bf x}^{r,i-1})\| = \|E_{\cdot i}\| |x^{r+1}_i - 
x^r_i| \geq \|E_{\cdot i}\| \epsilon, 
\end{equation}

\noindent
where $E_{\cdot i}$ denotes the $i^{th}$ column of $E$. Since $\{{\bf t}^{r,i}\}_{r \in \mathcal{R}}$ and 
$\{{\bf t}^{r,i-1}\}_{r \in \mathcal{R}}$ are bounded, we assume without loss of generality that there is 
a further subsequence $\mathcal{R}'$ of $\mathcal{R}$ such that $\{{\bf t}^{r,i}\}_{r \in \mathcal{R}'}$ 
and $\{{\bf t}^{r,i-1}\}_{r \in \mathcal{R}'}$ converge to ${\bf t}'$ and ${\bf t}''$ respectively. It follows 
by (\ref{eq7.1}) that ${\bf t}' \neq {\bf t}''$. Since ${\bf t}', {\bf t}'' \in C_g$, it follows by the continuity 
of $g$ that 
\begin{equation} \label{eq7}
\{g({\bf t}^{r,i})\}_{r \in \mathcal{R}'} \rightarrow g({\bf t}'), \; \{g({\bf t}^{r,i-1})\}_{r \in \mathcal{R}'} 
\rightarrow g({\bf t}''). 
\end{equation}

\noindent
It follows by the definition of $f$ that 
\begin{equation} \label{eq8}
\left\{\sum_{j \in S} |x^{r,i}_j|\right\}_{r \in \mathcal{R}'} \rightarrow f^\infty - g(t'), \; 
\left\{\sum_{j \in S} |x^{r,i-1}_j|\right\}_{r \in \mathcal{R}'} \rightarrow f^\infty - g(t''). 
\end{equation}

\noindent
Since ${\bf x}^{r,i}$ is obtained from ${\bf x}^{r, i-1}$ by minimizing along the $i^{th}$ coordinate, 
the convexity of $f$ yields 
\begin{eqnarray}
f({\bf x}^{r,i}) 
&\leq& f \left( \frac{{\bf x}^{r,i} + {\bf x}^{r,i-1}}{2} \right) = g \left( \frac{{\bf t}^{r,i} + 
{\bf t}^{r,i-1}}{2} \right) + \frac{\sum_{j \in S} |x^{r,i}_j + x^{r,i-1}_j|}{2} \nonumber\\
&\leq& g \left( \frac{{\bf t}^{r,i} + {\bf t}^{r,i-1}}{2} \right) + \frac{\sum_{j \in S} |x^{r,i}_j| + 
|x^{r,i-1}_j|}{2} \label{eq9}, 
\end{eqnarray}

\noindent
for every $r \in \mathcal{R}'$. Using the continuity of $g$, (\ref{eq8}) and passing to the limit as 
$r \rightarrow \infty, r \in \mathcal{R}'$, we obtain 
$$
f^\infty \leq f^\infty + g \left( \frac{{\bf t}' + {\bf t}''}{2} \right) - \frac{g({\bf t}') + g({\bf t}'')}{2}. 
$$

\noindent
The above yields a contradiction to the strict convexity of $g$ on $C_g$. \hfill$\Box$ 

\bigskip

\noindent
Using the result from Lemma \ref{lem4} above, we now proceed to prove that $\{{\bf t}^{r,i}\}_{r \geq 
0}$ converges to ${\bf t}^*$ for every $0 \leq i \leq n$, and then use this to establish that the 
sequence of differences between the successive iterates produced by the cyclic coordinatewise 
descent algorithm fo $f_1$ is 
square-summable. 
\begin{lemma} \label{lem5}
\begin{enumerate}[(a)]
\item For every $0 \leq i \leq n$, 
\begin{equation} \label{eq10}
\|{\bf t}^{r,i} - {\bf t}^*\| \rightarrow 0, 
\end{equation}

\noindent
as $r \rightarrow \infty$. 
\item 
$$
\sum_{r=0}^\infty \|{\bf x}^r - {\bf x}^{r+1}\|^2 < \infty. 
$$
\end{enumerate}
\end{lemma}

\noindent
{\it Proof} (a) Fix $i$ between $0$ to $n$ arbitrarily. Since $\{{\bf t}^{r,i}\}_{r \geq 0}$ is bounded, it 
has at least one limit point. Let ${\bf t}^\infty$ be an arbitrarily chosen limit point. Hence, there exists 
a subsequence $\mathcal{R}$ of $\mathbb{N}$ such that $\{{\bf t}^{r,i}\}_{r \in \mathcal{R}}$ converges 
to ${\bf t}^\infty$. Note that ${\bf t}^\infty \in C_g$. Hence, $g$ is continuously differentiable in an 
open set around ${\bf t}^\infty$. 

Note that for every $j \neq i$, 
$$
\|{\bf x}^{r,j} - {\bf x}^{r,i}\| = \sqrt{\sum_{k=\min(i,j)+1}^{\max(i,j)} |x^{r+1}_k - x^r_k|^2} \leq \|{\bf x}^{r+1} - {\bf x}^r\|. 
$$

\noindent
It follows by Lemma \ref{lem4} that $\|{\bf x}^{r,j} - {\bf x}^{r,i}\| \rightarrow 0$ as $r 
\rightarrow \infty$ for every $0 \leq j \leq n$. Hence, we have $\|{\bf t}^{r,j} - {\bf t}^{r,i}\| 
\rightarrow 0$ as $r \rightarrow \infty$ for every $0 \leq j \leq n$. It follows that 
\begin{equation} \label{eq10.1}
\{{\bf t}^{r,j}\}_{r \in \mathcal{R}} \rightarrow {\bf t}^\infty 
\end{equation}

\noindent
for every $0 \leq j \leq n$. Let $d^\infty = E^T \nabla g({\bf t}^\infty)$. It follows that 
\begin{equation} \label{eq11}
\{d({\bf x}^{r,j}\}_{r \in \mathcal{R}} = \{E^T \nabla g({\bf t}^{r,j})\}_{r \in \mathcal{R}} \rightarrow 
d^\infty 
\end{equation}

\noindent
as for every $0 \leq j \leq n$. By (\ref{eq4}) and (\ref{eq5}), it follows that for every $r \in 
\mathcal{R}$, 
\begin{equation} \label{eq12}
x^{r+1}_i  = x^{r,i}_i = sign(x^{r,i}_i - d_i ({\bf x}^{r,i})) \max(|x^{r,i}_i - d_i ({\bf x}^{r,i})| - \lambda, 0), 
\end{equation}

\noindent
for $i \in S$, and 
\begin{equation} \label{eq13}
x^{r+1}_i = x^{r,i}_i = \max(0, x^{r,i}_i - d_i ({\bf x}^{r,i})). 
\end{equation}

\noindent
for $i \in S^c$. By the arguments in the proof of Lemma \ref{lem2}, it follows that 
\begin{eqnarray}
& & d_i ({\bf x}^{r,i}) + \lambda sign(x^{r,i}_i) = 0 \; \mbox{ if } x_i \neq 0, i \in S, \label{eq13.1}\\
& & |d_i ({\bf x}^{r,i})| \leq \lambda \; \mbox{ if } x_i = 0, i \in S. \label{eq13.2} 
\end{eqnarray}

\noindent
It follows from (\ref{eq11}) and (\ref{eq12}) that $|d^\infty_i| \leq \lambda$ for $i \in S$. Since 
$x^{r,i}_i \geq 0$ for $i \in S^c$, it follows from (\ref{eq11}) and (\ref{eq13}) that $d^\infty_i \geq 0$ 
for $i \in S^c$. If $i \in S$ and $|d^\infty_i| < \lambda$, then $|d_i ({\bf x}^{r,i})| < \lambda$ for large 
enough $r$. It follows that 
\begin{equation} \label{eq14}
x^{r+1}_i = x^{r,i}_i = 0. 
\end{equation}

\noindent
If $i \in S^c$ and $d^\infty_i > 0$, then $d_i ({\bf x}^{r,i}) > 0$ for large enough $r$. It follows that 
\begin{equation} \label{eq15} 
x^{r+1}_i = x^{r,i}_i = 0. 
\end{equation}

\noindent
For each $r \in \mathcal{R}$, consider the 
linear system 
\begin{equation} \label{eq16}
E {\bf x} = {\bf t}^{r+1}, \; x_j = x^{r+1}_j \; \forall j \in S \mbox{ and } j \in S^c \mbox{ with } 
d^\infty_j > 0, \; {\bf x} \in \mathcal{X}. 
\end{equation}

\noindent
This a consistent system of equations since ${\bf x}^{r+1}$ is a solution. Fix any ${\bar {\bf x}} \in 
\mathcal{X}$. By Lemma \ref{lem3}, for every $r \in \mathcal{R}$, there exists a solution ${\bf y}^r$ of 
this linear system satisfying 
\begin{equation} \label{eq17}
\|{\bar {\bf x}} - {\bf y}^r\| \leq \theta \left( \|E {\bar {\bf x}} - {\bf t}^{r+1}\| + \sum_{j \in S} 
|\bar{x}_j - x^{r+1}_j| + \sum_{j \in S^c: d^\infty_j > 0} |\bar{x}_j - x^{r+1}_j| \right), 
\end{equation}

\noindent
where $\theta$ is a constant depending on $E$ only. Note that by (\ref{eq6}), $\{f_1 
({\bf x}^{r,i})\}_{r \in \mathbb{R}, 1 \leq i \leq n}$ is bounded above. By 
(\ref{eq6.2}), we get that $\{{\bf t}^{r,i}\}_{r \in \mathcal{R}, 1 \leq i \leq n}$ is 
bounded. Hence, $\{g({\bf t}^{r,i})\}_{r \in \mathcal{R}, 1 \leq i \leq n}$ is bounded 
below. It follows by the definition of $f$ that $\{\sum_{j \in S} |x^{r,i}_j|\}_{r 
\in \mathcal{R}, 1 \leq i \leq n}$ is bounded above. Hence, the right hand side of 
(\ref{eq17}) is bounded for all $r \in \mathcal{R}$. It follows by (\ref{eq10.1}), 
(\ref{eq14}) and (\ref{eq15}) that $\{{\bf y}^r\}_{r \in \mathcal{R}}$ is bounded, and 
that every limit point of $\{{\bf y}^r\}_{r \in \mathcal{R}}$, say ${\bf y}^\infty$, 
satisfies 
\begin{equation} \label{eq17.1}
E {\bf y}^\infty = {\bf t}^\infty, \; y^\infty_j = 0 \; \forall j \in S \mbox{ with } |d^\infty_j| < \lambda  
\mbox{ and } j \in S^c \mbox{ with } d^\infty_j > 0, \; {\bf y}^\infty \in \mathcal{X}. 
\end{equation}

\noindent
Since $E {\bf y}^\infty = {\bf t}^\infty$, we obtain $d({\bf y}^\infty) = E^T \nabla g({\bf t}^\infty) = 
d^\infty$. 

Note that if $j \in S$ and $d^\infty_j = \lambda$, it follows from (\ref{eq11}), (\ref{eq13.1}) 
and (\ref{eq13.2}) that $x^{r+1}_j = x^{r,j}_j \leq 0$ for large enough $r$. Since ${\bf y}^r$ satisfies 
(\ref{eq16}), and ${\bf y}^\infty$ is a limit point, it follows that $y^\infty_j \leq 0$. Hence, 
\begin{equation} \label{eq18}
y^\infty_j = sign(y^\infty_j - d_j ({\bf y}^\infty)) \max(|y^\infty_j - d_j ({\bf y}^\infty)| - \lambda, 0). 
\end{equation}

\noindent
If $j \in S$ and $d^\infty_j = -\lambda$, a similar argument as above implies that $y^\infty_j$ 
satisfies (\ref{eq18}). If $j \in S^c$ and $d^\infty_j = d_j ({\bf y}^\infty) > 0$, then it follows by 
(\ref{eq15}) that ${\bf y}^infty_j = 0$. Hence, 
\begin{equation} \label{eq19}
y^\infty_j = \max(y^\infty_j - d_j ({\bf y}^\infty), 0). 
\end{equation}

\noindent
If $j \in S^c$ and $d^\infty_j > 0$, then (\ref{eq19}) holds trivially. It follows from Lemma \ref{lem2}, 
(\ref{eq17.1}), (\ref{eq18}) and (\ref{eq19}) that ${\bf y}^\infty \in \mathcal{X}^*$. It follows by 
Lemma \ref{lem1} that ${\bf t}^\infty = E {\bf y}^\infty = {\bf t}^*$. Since ${\bf t}^\infty$ is an 
arbitrarily chosen limit point of $\{{\bf t}^{r,i}\}_{r \geq 0}$, it follows that for every $0 \leq i \leq n$, 
$$
\|{\bf t}^{r,i} - {\bf t}^*\| \rightarrow 0 
$$

\noindent
as $r \rightarrow \infty$. This establishes part (a). 

\medskip

\noindent
(b) By Lemma \ref{lem5}, it follows that for $r$ sufficiently large, $E {\bf x}^{r,i} \in 
\mathcal{U}^*$ for every $0 \leq i \leq n$. Consider any such $r$. For every $1 \leq i \leq n$, a 
second order Taylor series expansion along the $i^{th}$ coordinate leads to 
$$
g(E {\bf x}^{r,i-1}) - g(E {\bf x}^{r,i}) = \nabla g(E {\bf x}^{r,i})^T E_{\cdot i} (x^{r,i-1}_i - x^{r,i}_i) + 
E_{\cdot i}^T \nabla^2 g(E \tilde{\bf x}^{r,i}) E_{\cdot i} (x^{r,i-1}_i - x^{r,i}_i)^2, 
$$

\noindent
where $\tilde{\bf x}^{r,i}$ is a convex combination of ${\bf x}^{r,i-1}$ and ${\bf x}^{r,i}$. 
Since $U^*$ is an open ball containing ${\bf t}^{r,i-1} = E {\bf x}^{r,i-1}$ and ${\bf t}^{r,i} = E 
{\bf x}^{r,i}$, we conclude that $E \tilde{\bf x}^{r,i}$ is contained in $U^*$. Note that $d_i 
({\bf x}^{r,i}) = \nabla g(E {\bf x}^{r,i})^T E_{\cdot i}$. It follows from (\ref{eq1}) that 
\begin{equation} \label{eq20}
g(E {\bf x}^{r,i-1}) - g(E {\bf x}^{r,i}) \geq d_i ({\bf x}^{r,i}) (x^{r,i-1}_i - x^{r,i}_i) + 
\sigma \|E_{\cdot i}\|^2 (x^{r,i-1}_i - x^{r,i}_i)^2. 
\end{equation}

\noindent
for every $1 \leq i \leq n$. 

Fix $i \in S$ arbitrarily. Let $h(x) := |x|$. For any subderivative $\delta$ of the function $h$ at 
$x^{r,i}_i$, we have 
\begin{equation} \label{eq21}
|x^{r,i-1}_i| - |x^{r,i}_i| \geq \delta (x^{r,i-1}_i - x^{r,i}_i). 
\end{equation}

\noindent
Note that ${\bf x}^{r,i-1}$ and ${\bf x}^{r,i}$ only differ in the $i^{th}$ coordinate. Using the 
definition of $f$ along with (\ref{eq20}) and (\ref{eq21}), it follows that 
$$
f({\bf x}^{r,i-1}) - f({\bf x}^{r,i}) \geq (d_i ({\bf x}^{r,i}) + \lambda \delta) (x^{r,i-1}_i - x^{r,i}_i) + 
\sigma \|E_{\cdot i}\|^2 (x^{r,i-1}_i - x^{r,i}_i)^2. 
$$

\noindent
Note that if $x^{r_i}_i \neq 0$, then $\delta = sign(x^{r,i}_i)$. If $x^{r,i}_i = 0$, then any $\delta 
\in [-1,1]$ is a valid subderivative choice for $h$, and ${\bf x}^{r,i}$ is obtained from 
${\bf x}^{r,i-1}$ by minimizing $f$ along the $i^{th}$ coordinate. Using these observations, we 
conclude that it is always possible to choose $\delta$ such that $d_i ({\bf x}^{r,i}) + \lambda \delta 
= 0$. Hence, for every $i \in S$ 
\begin{equation} \label{eq22}
f({\bf x}^{r,i-1}) - f({\bf x}^{r,i}) \geq \sigma \left( \min_{1 \leq j \leq n} \|E_{\cdot j}\|^2 \right) 
(x^{r,i-1}_i - x^{r,i}_i)^2. 
\end{equation}

\noindent
Fix $i \in S^c$ arbitrarily. By (\ref{eq5}), $d_i ({\bf x}^{r,i}) \geq 0$. Suppose $d_i ({\bf x}^{r,i}) = 0$. 
Since ${\bf x}^{r,i-1}$ and ${\bf x}^{r,i}$ only differ in the $i^{th}$ coordinate, it follows by 
(\ref{eq20}), the definition of $f$ that (\ref{eq22}) holds in this case. Suppose $d_i ({\bf x}^{r,i}) > 
0$. By (\ref{eq5}), $x^{r,i}_i = 0$. Since $x^{r,i-1}_i \geq 0$, it follows that $d_i ({\bf x}^{r,i}) 
(x^{r,i-1}_i - x^{r,i}_i) \geq 0$. Hence, (\ref{eq22}) holds in this case. 

Adding (\ref{eq22}) over $i = 1,2, \cdots, n$, we obtain 
$$
 f({\bf x}^r) - f({\bf x}^{r+1}) \geq \sigma \left( \min_{1 \leq j \leq n} \|E_{\cdot j}\|^2 \right) 
\sum_{i=1}^n (x^{r,i-1}_i - x^{r,i}_i)^2 = \sigma \left( \min_{1 \leq j \leq n} \|E_{\cdot j}\|^2 
\right) \|{\bf x}^r - {\bf x}^{r+1}\|^2. 
$$

\noindent
The result follows by noting that $f({\bf x}^r) \downarrow f^\infty > -\infty$ as $r \rightarrow 
\infty$ and that $\min_{1 \leq j \leq n} \|E_{\cdot j}\|^2 > 0$ as $E$ has no zero column. 
\hfill$\Box$ 

\bigskip

\noindent
Although the square-summability, established above, is an important step towards 
proving convergence, further arguments are needed to establish convergence of the sequence 
of iterates generated by Algorithm 1. It follows by Lemma \ref{lem5} and the continuity of 
$\nabla g$ at ${\bf t}^*$ that 
\begin{equation} \label{eq23}
{\bf d} ({\bf x}^{r,i}) \rightarrow {\bf d}^* 
\end{equation}

\noindent
as $r \rightarrow \infty$ for every $1 \leq i \leq n$.  The next lemma establishes 
that for each $i$, $x^r_i$ has the same sign for sufficiently large $r$. 
\begin{lemma} \label{lem8}
\begin{enumerate}[(a)]
\item For all $r$ sufficiently large, $x_i^r = 0$ for all $i \in S$ with $|d_i^*| < \lambda$ and 
for all $i \in S^c$ with $d_i^* > 0$. 
\item For all $r$ sufficiently large, $x_i^r \leq 0$ for all $i \in S$ with $d_i^* = \lambda$, 
and $x_i^r \geq 0$ for all $i \in S$ with $d_i^* = -\lambda$. 
\end{enumerate}
\end{lemma}

\noindent
{\it Proof} The proof of part (a) follows by using exactly the same arguments as those 
leading to (\ref{eq14}) and (\ref{eq15}). We now prove part (b). Let $i \in S$ with $d_i^* = 
\lambda$. Note that by (\ref{eq4}), $x^r_i = sign(x^r_i - d_i ({\bf x}^{r-1,i})) \max(|x^r_i 
- d_i ({\bf x}^{r-1,i})| - \lambda, 0)$. By exactly the same arguments as in the proof of 
Lemma \ref{lem2}, it follows that $d_i ({\bf x}^{r-1,i}) + \lambda sign(x^r_i) = 0$ if 
$x_i^r \neq 0$. By (\ref{eq23}), it follows that $d_i ({\bf x}^{r-1,i}) \rightarrow \lambda$ as $r 
\rightarrow \infty$. Hence, for sufficiently large r, $x_i^r \neq 0$ implies that $x_i^r > 0$. The 
other case follows similarly. \hfill$\Box$ 

\medskip

\noindent
For every ${\bf x} \in \mathcal{X}$, define the function $\phi$ as follows: 
$$
\phi ({\bf x}) = \min_{{\bf x}^* \in \mathcal{X}^*} \|{\bf x} - {\bf x}^*\|. 
$$

\noindent
Hence, $\phi({\bf x})$ is the distance of ${\bf x}$ from the closed convex set $\mathcal{X}^*$. 
The goal of the next lemma is to establish that the sequence of iterates 
$\{{\bf x}^r\}_{r \geq 0}$ approaches the boundary of $\mathcal{X}^*$. Although,this lemma 
is a useful component of the convergence proof of the CCM algorithm for $f_1$, it clearly is not 
sufficient to establish convergence. For example, consider a sequence which alternatively takes two distinct values at the 
boundary of a set. It easily follows that the distance of the sequence from the boundary of that set is always 
zero, but the sequence still does not converge. 
\begin{lemma} \label{lem11}
\begin{enumerate}[(a)]
\item Let $\lambda > 0$. Then 
\begin{equation} \label{eq24}
|sign(a) \max(|a|-\lambda, 0) - sign(b) \max(|b|-\lambda, 0)|  \leq |a - b|, 
\end{equation}

\noindent
for all $a,b \in \mathbb{R}$. 
\item If $i \in S$, then 
$$
x^r_i - sign(x^r_i - d_i ({\bf x}^r)) \max(|x^r_i - d_i ({\bf x}^r)| - \lambda, 0) \rightarrow 0 
$$

\noindent
as $r \rightarrow \infty$. If $i \in S^c$, then 
$$
x^r_i - \max(x^r_i - d_i ({\bf x}^r), 0) \rightarrow 0 
$$

\noindent
as $r \rightarrow \infty$. 
\item $$
\phi({\bf x}^r) \rightarrow 0 \mbox{ as } r \rightarrow \infty. 
$$
\end{enumerate}
\end{lemma}

\noindent
{\it Proof} (a) We consider various cases:\\

\noindent
Suppose $|a| \leq \lambda$ and $|b| \leq \lambda$. Then the $\max(|a| - \lambda, 0) = 
\max(|b| - \lambda, 0) = 0$. Hence, (\ref{eq24}) holds. Suppose $|a| \leq \lambda$ and 
$|b| > \lambda$. Then 
\begin{eqnarray*}
|sign(a) \max(|a|-\lambda, 0) - sign(b) \max(|b|-\lambda, 0)| 
&=& |b| - \lambda\\
&\leq& |b| - |a|\\
&\leq& |b-a|. 
\end{eqnarray*}

\noindent
The last step follows by the triangle inequality. The case $|a| > \lambda$ and $|b| \leq \lambda$ 
can be analyzed similarly. Next, suppose that $|a| > \lambda$, $|b| > \lambda$ and $sign(a) = 
sign(b)$. Then 
\begin{eqnarray*}
|sign(a) \max(|a|-\lambda, 0) - sign(b) \max(|b|-\lambda, 0)| 
&=& |a - b + \lambda sign(b) - \lambda sign (a)|\\
&=& |a - b|. 
\end{eqnarray*}

Finally, we consider the case when $|a| > \lambda$, $|b| > \lambda$ and $sign(a) \neq sign(b)$. 
Without loss of generality, let $a < 0$ and $b > 0$. Then 
\begin{eqnarray*}
|sign(a) \max(|a|-\lambda, 0) - sign(b) \max(|b|-\lambda, 0)| 
&=& |a - b + 2 \lambda|\\
&=& b + |a| - 2 \lambda\\
&<& b + |a|\\
&=& |a - b|. 
\end{eqnarray*}
 
\medskip

\noindent
(b) Note that by (\ref{eq4}), $x^r_i = sign(x^r_i - d_i ({\bf x}^{r-1,i})) \max(|x^r_i - d_i 
({\bf x}^{r-1,i})| - \lambda, 0)$. It follows by Lemma \ref{lem8}, part (a) of this lemma, and 
(\ref{eq23}) that 
$$
|x^r_i - sign(x^r_i - d_i ({\bf x}^r)) \max(|x^r_i - d_i ({\bf x}^r)| - \lambda, 0)| \leq |d_i 
({\bf x}^{r-1,i}) - d_i ({\bf x}^r)| \rightarrow 0 
$$

\noindent
as $r \rightarrow \infty$. Similarly, by (\ref{eq5}), $x^r_i = \max(x^r_i - d_i ({\bf x}^{r-1,i}), 0)$. 
Note that for any $a,b \in \mathbb{R}$, $|\max(a,0) - \max(b,0)| < |a - b|$. It follows by 
(\ref{eq23}) that 
$$
|x^r_i - \max(x^r_i - d_i ({\bf x}^r), 0)| \leq |d_i ({\bf x}^{r-1,i}) - d_i ({\bf x}^r)| \rightarrow 0. 
$$

\noindent
as $r \rightarrow \infty$. 

\medskip

\noindent
(c) By Lemma \ref{lem1}, (\ref{kkt1}), (\ref{kkt2}) and (\ref{kkt3}), it follows that 
$\mathcal{X}^*$ is the solution set of the linear system of equations given by 
\begin{eqnarray*}
& & E {\bf y} = {\bf t}^*, {\bf y} \in \mathcal{X},\\
& & y_i = 0 \mbox{ if } i \in S, |d_i^*| < \lambda,\\
& & y_i \leq 0 \mbox{ if } i \in S, d_i^* = \lambda,\\
& & y_i \geq 0 \mbox{ if } i \in S, d_i^* = -\lambda,\\
& & y_i = 0 \mbox{ if } i \in S^c, d_i^* > 0. 
\end{eqnarray*}

\noindent
Since $\mathcal{X}^*$ is non-empty, by Lemma \ref{lem3} and Lemma \ref{lem8}, for sufficiently 
large $r$, there exists ${\bf y}^r \in \mathcal{X}^*$ such that 
\begin{eqnarray} 
\|{\bf x}^r - {\bf y}^r\| 
&\leq& \theta \left( \|E {\bf x}^r - {\bf t}^*\| + \sum_{i \in S, d_i^* = \lambda} (x^r_i)^+ + 
\sum_{i \in S, d_i^* = -\lambda} (-x^r_i)^+  + \sum_{i \in S, |d_i^*| < \lambda} |x^r_i| \right) + 
\nonumber\\
& & \theta \sum_{i \in S^c, d_i^* > 0} |x^r_i| \nonumber\\
&=& \theta \|E {\bf x}^r - {\bf t}^*\| \label{eq25}
\end{eqnarray}

\noindent
where $\theta$ is a constant only depending on $E$. The result follows by the definition of $\phi$ and 
Lemma \ref{lem5}. \hfill$\Box$ 

\bigskip

\noindent
Let 
\begin{eqnarray*}
& & I_1^* := \{i \in S: d_i^* = \lambda\},\\
& & I_2^* := \{i \in S: d_i^* = -\lambda\},\\
& & I_3^* := \{i \in S: |d_i^*| < \lambda\},\\
& & I_4^* := \{i \in S^c: d_i^* = 0\},\\
& & I_5^* := \{i \in S^c: d_i^* > 0\}. 
\end{eqnarray*}

\noindent
For ${\bf x} \in \mathbb{R}^n$ and $M \subseteq \{1,2, \cdots, n\}$, let ${\bf x}_M := (x_i)_{i \in M}$. 
By Lemma \ref{lem8} and (\ref{eq5}) that there exists an $r_0 > 0$ such that 
\begin{equation}
{\bf x}^r_{I_3^* \cup I_5^*} = {\bf 0}, \; ({\bf x}^r_{I_1^*})^+ = {\bf 0}, \; (-{\bf x}^r_{I_2^* \cup 
I_4^*})^+ = {\bf 0} \label{eq28} 
\end{equation}

\noindent
for every $r \geq r_0$. The following lemma provides a crucial identity which will play 
in important role in the last leg of the convergence proof. 
\begin{lemma} \label{lem13}
There exists $\omega > 0$ such that 
$$
\|E {\bf x}^r - {\bf t}^*\| \leq \omega \|{\bf x}^r - {\bf x}^{r+1}\| 
$$

\noindent
for every $r \geq r_0$. 
\end{lemma}

\noindent
{\it Proof} Consider arbitrary (possibly empty) subsets $I_1, I_2, I_3$ of $S$ and $I_4, I_5$ of 
$S^c$, and let $\mathcal{R}$ denote the set of indices $r \geq r_0$ for which 
\begin{eqnarray}
d_i ({\bf x}^{r,i}) = \lambda \; \forall i \in I_1, \label{eq29}\\
d_i ({\bf x}^{r,i}) = -\lambda \; \forall i \in I_2, \label{eq30}\\
|d_i ({\bf x}^{r,i})| < \lambda \; \forall i \in I_3, \label{eq31}\\
d_i ({\bf x}^{r,i}) = 0 \; \forall i \in I_4, \label{eq32}\\
d_i ({\bf x}^{r,i}) > 0 \; \forall i \in I_5. \label{eq33}
\end{eqnarray}

\noindent
Note that $|d_i ({\bf x}^{r,i}| \leq \lambda$ for $i \in S$, and $d_i ({\bf x}^{r,i}) \geq 0$ for $i \in 
S^c$. Hence, $I_1 \cup I_2 \cup I_3 = S$ and $I_4 \cup I_5 = S^c$. Suppose we are able to show 
that there exists a constant $\omega_{I_1, I_2, I_3, I_4, I_5} > 0$ such that 
\begin{equation} \label{eq33.1}
\|E {\bf x}^r - {\bf t}^*\| \leq \omega_{I_1, I_2, I_3, I_4, I_5} \|{\bf x}^r - {\bf x}^{r+1}\|, 
\end{equation}

\noindent
for every $r \in \mathcal{R}$. Since every $r \geq r_0$ belongs to $\mathcal{R}$ corresponding to 
some choice of $\{I_j\}^{1 \leq j \leq 5}$, and the number of distinct choices of $\{I_j\}_{1 \leq j 
\leq 5}$ is finite, it would immediately imply that the lemma holds with 
$$
\omega = \max_{I_1, I_2, I_3, I_4, I_5} \omega_{I_1, I_2, I_3, I_4, I_5}. 
$$

\noindent
Hence, we now establish (\ref{eq33.1}). Note that if ${\bf x}^{r+1} = {\bf x}^r$, then by 
Lemma \ref{lem1}, Lemma \ref{lem2}, ({\ref{eq4}) and (\ref{eq5}), it follows that ${\bf x}^r \in 
\mathcal{X}^*$ and $E {\bf x}^r = {\bf t}^*$. Hence, if $\mathcal{R}$ is empty or finite, then 
the result holds trivially. Hence, we assume that $\mathcal{R}$ is infinite. It follows by 
(\ref{eq4}), (\ref{eq5}) and (\ref{eq29})-(\ref{eq33}) that 
\begin{equation} \label{eq34}
{\bf x}^{r+1}_{I_3 \cup I_5} = {\bf 0}, \;  ({\bf x}^{r+1}_{I_1})^+ =  {\bf 0}, \; (-{\bf x}^{r+1}_{I_2 \cup 
I_4})^+ = {\bf 0}. 
\end{equation}

\noindent
Consider the linear system 
\begin{equation} \label{eq35}
{\bf y}_{I_3 \cup I_5} = {\bf 0}, {\bf y} \in \mathcal{X}^*. 
\end{equation}

\noindent
By an argument very similar to the one following \cite[eq. (B.7)]{Luo:Tseng:1989}, it follows that 
the above linear system is consistent (essentially by noting that $\{{\bf y} \in \mathcal{X}: 
{\bf y}_{I_3 \cup I_5} = {\bf 0}\}$ and $\mathcal{X}^*$ are polyhedral sets, and proving that they 
get arbitrarily close to each other). It follows by Lemma \ref{lem1}, (\ref{kkt1}), (\ref{kkt2}) and 
(\ref{kkt3}) that the solution set of the linear system in (\ref{eq35}) is identical to the solution set 
of the following linear system. 
\begin{eqnarray*}
& & {\bf y}_{I_3 \cup I_5} = {\bf 0},\\
& & E {\bf y} = {\bf t}^*, {\bf y} \in \mathcal{X},\\
& & {\bf y}_{I_3^* \cup I_5^*} = {\bf 0},\\
& & ({\bf y}_{I_1^*})^+ = {\bf 0},\\
& & (-{\bf y}_{I_2^* \cup I_4^*})^+ = {\bf 0}. 
\end{eqnarray*}

\noindent
It follows by Lemma \ref{lem3} that, for every $r \in \mathcal{R}$, there exists a solution ${\bf y}^r$ 
to the above linear system satisfying 
\begin{equation} \label{eq36}
\|{\bf x}^r - {\bf y}^r\| \leq \kappa_1 (\|E {\bf x}^r - {\bf t}^*\| + \|{\bf x}^r_{I_3 \cup I_3^* \cup 
I_5 \cup I_5^*}\| + \|({\bf x}^r_{I_1^*})^+\| + \|(-{\bf x}^r_{I_2^* \cup I_4^*})^+\|), 
\end{equation}

\noindent
where $\kappa_1$ depends only on $E$. It follows by (\ref{eq28}) and (\ref{eq34}) that 
\begin{eqnarray}
\|{\bf x}^r - {\bf y}^r\| 
&\leq& \kappa_1 (\|E {\bf x}^r - {\bf t}^*\| + \|{\bf x}^r_{I_3 \cup I_3^* \cup I_5 \cup I_5^*} - 
{\bf x}^{r+1}_{I_3 \cup I_3^* \cup I_5 \cup I_5^*}\|) \nonumber\\
&\leq& \kappa_1 (\|E {\bf x}^r - {\bf t}^*\| + \|{\bf x}^r - {\bf x}^{r+1}\|). \label{eq37}
\end{eqnarray}

\noindent
For any $m \times n$ matrix $A$ and $M \subseteq \{1,2, \cdots, n\}$, let $A_M := ((A_{ij}))_{i 
\in \{1 \leq i \leq m, j \in M}$. Let $I = I_1 \cup I_2 \cup I_4$. Note that by (\ref{eq34}) 
$\|E_{I^c} ({\bf x}^r_{I^c} - {\bf y}^r_{I^c})\| = \|E_{I^c} ({\bf x}^r_{I^c} - {\bf x}^{r+1}_{I^c})\| 
\leq \|E\| \|{\bf x}^r - {\bf x}^{r+1}\|$. It follows by Lemma \ref{lem1} and (\ref{eq37}) that 
\begin{eqnarray}
\|{\bf x}^r_{I^c} - {\bf y}^r_{I^c}\| 
&\leq& \|{\bf x}^r - {\bf y}^r\| \nonumber\\
&\leq& \kappa_1 (\|E ({\bf x}^r - {\bf y}^r)\| + \|{\bf x}^r - {\bf x}^{r+1}\|) \nonumber\\
&\leq& \kappa_1 ((1 + \|E\|) \|{\bf x}^r - {\bf x}^{r+1}\| + \|E_I ({\bf x}^r_I - {\bf y}^r_I)\|) 
\nonumber\\
&\leq& \kappa_1 (1 + \|E\|) (\|{\bf x}^r - {\bf x}^{r+1}\| + \|E_I ({\bf x}^r_I - {\bf y}^r_I)\|) 
\label{eq38}
\end{eqnarray}

\noindent
Let ${\bf c} \in \mathbb{R}^n$ be such that ${\bf c}_{I_1} = \lambda$, ${\bf c}_{I_2} = 
-\lambda$, and all the other entries of ${\bf c}$ are equal to zero. It follows by 
Lemma \ref{lem1}, (\ref{eq29}), (\ref{eq30}) and (\ref{eq32}) that 
$$
\|{\bf d}_I ({\bf x}^r) - {\bf c}_I\| = \|{\bf d}_I ({\bf x}^r) - {\bf d}_I^*\| = \|{\bf d}_I ({\bf x}^r) - 
{\bf d}_I ({\bf y}^r)\| = \|(E_I)^T \nabla g(E {\bf x}^r) - (E_I)^T \nabla g(E {\bf y}^r)\|, 
$$

\noindent
and 
$$
|d_i ({\bf x}^r) - c_i| = |d_i ({\bf x}^r) - d_i ({\bf x}^{r,i})| = |E_{\cdot i}^T \nabla g(E {\bf x}^{r}) - 
E_{\cdot i}^T \nabla g(E {\bf x}^{r,i})|. 
$$

\noindent
for every $i \in I$. The result now follows by exactly the same arguments as in 
\cite{Luo:Tseng:1989} (from \cite[eq. (B.9)]{Luo:Tseng:1989} to the end of the proof of 
\cite[Lemma B.3]{Luo:Tseng:1989}, using $\|{\bf d}_I ({\bf x}^r) - {\bf c}_I\|$ in place of 
$\|{\bf d}_I ({\bf x}^r)\|$, and replacing $E$ by $E^T$ throughout). \hfill$\Box$ 

\bigskip

\noindent
We now invoke two matrix-theoretic results from Luo and Tseng \cite{Luo:Tseng:1989}. Let $M = 
E^T \nabla^2 g({\bf t}^*) E$. By (A4) and the assumption that $E$ has no zero column, it follows 
that $m_{ii} > 0$ for every $1 \leq i \leq n$. For any $J, \tilde{J} \subseteq \{1,2, \cdots, n\}$, 
let $M_{J \tilde{J}} := (M_{ij})_{i \in J, j \in \tilde{J}}$, and $|J|$ denote the cardinality of $J$. The 
following lemma is provided in Luo and Tseng \cite{Luo:Tseng:1989}, and exploits the fact 
that $M$ is symmetric positive semi-definite. 
\begin{lemma}[Luo and Tseng \cite{Luo:Tseng:1989}] \label{lem14}
Let $J \subseteq \{1,2, \cdots, n\}$. Then $Span(M_{JJ^c}) \subseteq Span(M_{JJ})$. 
\end{lemma}


\noindent
Let $B$ denote the lower triangular portion of $M$, and $C = M - B$ denote the strictly upper 
triangular portion of $M$ (hence the diagonal entries of $C$ are zero). We use 
the following lemma from \cite{Luo:Tseng:1989}. 
\begin{lemma}[Luo and Tseng \cite{Luo:Tseng:1989}] \label{lem15}
\begin{enumerate}[(a)]
\item For any nonempty $J \subseteq \{1,2, \cdots, n\}$, there exist $\rho_J \in (0,1)$ and 
$\tau_J > 0$ such that 
$$
\left\| (I - M_{JJ} (B_{JJ})^{-1})^k {\bf z} \right\| \leq \tau_J (\rho_J)^k \|{\bf z}\|, \; \forall k \geq 
1, \; \forall {\bf z} \in Span(M_{JJ}). 
$$
\item There exists a $\Delta \geq 1$ such that, for any nonempty $J \subseteq \{1,2, \cdots, 
n\}$, 
$$
\left\| (I - M_{JJ} (B_{JJ})^{-1})^k {\bf z} \right\| \leq \Delta \|{\bf z}\|, \; \forall k \geq 1, \; \forall 
{\bf z} \in \mathbb{R}^{|J|}. 
$$
\end{enumerate}
\end{lemma}

\noindent
Let $I^* = I_1^* \cup I_2^* \cup I_4^*$, and 
$$
\beta = \max_{J \subseteq I^*} \sqrt{|J^c|} \left\{ \left( \frac{\tau_J \|(B_{JJ}^{-1}\| \|M_{JJ}\|}{1 - 
\rho_J} + \Delta + 1 \right) \|(B_{JJ})^{-1} B_{JJ^c}\| + \frac{\tau_J \|(B_{JJ})^{-1}\| \|M_{JJ}\|}{1 - 
\rho_J} \right\}. 
$$

\noindent
Recall tha by Lemma \ref{lem8} and (\ref{eq5}), there exists an $r_0 > 0$ such that 
$$
{\bf x}^r_{I_3^* \cup I_5^*} = {\bf 0}, \; ({\bf x}^r_{I_1^*})^+ = {\bf 0}, \; (-{\bf x}^r_{I_2^* \cup 
I_4^*})^+ = {\bf 0} 
$$

\noindent
for every $r \geq r_0$. For ${\bf x} \in \mathbb{R}^n$, let $\|{\bf x}\|_\infty = \max_{1 \leq i 
\leq n} |x_i|$. The next lemma is analogous to Lemma 9 of \cite{Luo:Tseng:1989}, and shows 
that the coordinates of ${\bf x}^r$ that stay away from zero, are influenced by the coordinates which 
eventually become zero only through the distance of these coordinates from zero. 
\begin{lemma} \label{lem16}
Consider any $J \subseteq I^*$. If for some two integers $s \geq t \geq r_0$ we have $x_i^r 
\neq 0$ for every $t+1 \leq r \leq s$ and $i \in J$, then, for any ${\bf x}^* \in \mathcal{X}^*$, 
there holds 
$$
\|{\bf x}_J^s - {\bf x}_J^*\| \leq \Delta \|{\bf x}_J^t - {\bf x}_J^*\| + \beta \max_{t \leq r \leq s} 
\|{\bf x}_{J^c}^r - {\bf x}_{J^c}^*\|_\infty + \mu \sum_{r=t}^{s-1} \|{\bf x}^r - 
{\bf x}^{r+1}\|^2, 
$$

\noindent
where $\mu$ is some positive constant which is independent of $s$ and $t$. 
\end{lemma}

\noindent
The proof of the lemma above is provided in the appendix. Let $\sigma_0 := 1$ and 
$$
\sigma_k = \Delta + 3+ \beta + (\beta + 1) \sigma_{k-1} + \mu, \; k = 1,2, \cdots, n. 
$$

\noindent
It follows from the above definition that $\sigma_k \geq 1$ for every $1 \leq k \leq n$, and is 
monotonically increasing with $k$. 

Fix $\delta > 0$ arbitrarily. By Lemma \ref{lem4}, Lemma \ref{lem5} and Lemma \ref{lem11}, 
there exists $r_1 > 0$ such that 
\begin{eqnarray}
\phi ({\bf x}^r) \leq \delta, \label{eq39}\\
\|{\bf x}^{r+1} - {\bf x}^r\| \leq \delta, \label{eq40}\\
\sum_{k=r}^\infty \|{\bf x}^k - {\bf x}^{k+1}\|^2 \leq \delta, \label{eq41} 
\end{eqnarray}

\noindent
for every $r \geq r_1$. 


The next three lemmas are analogous to \cite[Lemma 10]{Luo:Tseng:1:1989}, 
\cite[Lemma 11]{Luo:Tseng:1:1989}, and \cite[Lemma 9]{Luo:Tseng:1:1989} respectively. 
The crucial difference is that we consider absolute values of appropriate vector entries (as 
opposed to the lemmas in \cite{Luo:Tseng:1:1989}, which use the vector entries 
themselves). Recall that $I^* = I_1^* \cup I_2^* \cup I_4^*$. The proofs of all three 
lemmas are provided in the appendix. 
\begin{lemma} \label{lem17}
Fix $k \in \{1,2, \cdots, n\}$ arbitrarily. If for some nonempty $J \subset I^*$, and some 
intergers $t' > t \geq \max (r_0, r_1)$, we have 
\begin{eqnarray}
& & |x_i^t| > \sigma_k \delta, \; \forall i \in J, \label{eq42}\\
& & |x_i^r| \leq \sigma_{k-1} \delta, \; \forall i \notin J, \forall r = t, t+1, \cdots, t' - 1 
\label{eq43}, 
\end{eqnarray}

\noindent
then the following hold: 
\begin{enumerate}[(a)]
\item $|x_i^{t'}| > \sigma_{k-1} \delta$ for every $i \in J$. 
\item There exists an ${\bf x}^* \in \mathcal{X}^*$ such that 
$$
\|{\bf x}^r - {\bf x}^*\|_\infty \leq \sigma_k \delta, \; \forall r = t, t+1, \cdots, t'-1. 
$$
\end{enumerate}
\end{lemma}

\noindent
The next lemma extends the previous lemma by removing the assumption that the coordinates 
that start near zero remain near zero. 
\begin{lemma} \label{lem18}
Fix $k \in \{1,2, \cdots, n\}$ arbitrarily. If for some $J \subseteq I^*$ with $|J| \geq |I^*| - 
k + 1$ and some interger $t > \max(r_0, r_1)$ we have 
\begin{eqnarray}
& & |x_i^t| > \sigma_k \delta, \; \forall i \in J, \label{eq48}\\
& & |x_i^t| \leq \sigma_{k-1} \delta, \; \forall i \notin J, \label{eq49} 
\end{eqnarray}

\noindent
then there exists an ${\bf x}^* \in \mathcal{X}^*$ and a $\bar{t} \geq t$ satisfying 
\begin{equation} \label{eq50}
\|{\bf x}^r - {\bf x}^*\|_\infty \leq \sigma_k \delta, 
\end{equation}

\noindent
for every $r \geq \bar{t}$. 
\end{lemma}

\noindent
We use Lemma \ref{lem18} to establish the final lemma in our analysis. 
\begin{lemma} \label{lem19}
For any $\delta > 0$, there exists an ${\bf x}^* \in \mathcal{X}^*$ and $\hat{r} > 0$ such 
that 
\begin{equation} \label{eq53}
\|{\bf x}^r - {\bf x}^*\|_\infty \leq \sigma_n \delta + \delta, 
\end{equation}

\noindent
for every $r \geq \hat{r}$. 
\end{lemma}

\noindent
Using Lemma \ref{lem19}, we are now able to complete the proof of our meta-theorem, 
Theorem \ref{thm1}. 

\medskip

\noindent
{\it Proof of Theorem \ref{thm1}} Fix $\epsilon > 0$ arbitrarily. By Lemma \ref{lem19}, 
there exists ${\bf x}^* \in \mathcal{X}^*$ and $\hat{r} > 0$ such that 
$$
\|{\bf x}^r - {\bf x}^*\|_\infty < \frac{\epsilon}{2}, 
$$

\noindent
for every $r \geq \hat{r}$. Hence, for every $r_1, r_2 > \hat{r}$, we obtain by the 
triangle inequality that 
\begin{eqnarray*}
\|{\bf x}^{r_1} - {\bf x}^{r_2}\|_\infty 
&\leq& \|{\bf x}^{r_1} - {\bf x}^*\|_\infty + \|{\bf x}^{r_2} - {\bf x}^*\|_\infty\\
&<& \epsilon. 
\end{eqnarray*}

\noindent
It follows that the sequence of iterates $\{{\bf x}^r\}_{r \geq 0}$ form a Cauchy sequence. By 
Lemma \ref{lem11}, we conclude that $\{{\bf x}^r\}_{r \geq 0}$ converges to an element of 
$\mathcal{X}^*$. \hfill$\Box$ 
\newtheorem{remark}{Remark}
\begin{remark} \label{rem1}
Note that Theorem \ref{thm1} holds for any $m \times n$ matrix $E$ with non-zero columns, 
and any subset $S$ of $\{1,2, \cdots, n\}$. It follows that Theorem \ref{thm1} holds for an 
arbitrary permutation of the order in which the $n$ coordinates are updated in the 
cyclic coordinatewise descent algorithm. 
\end{remark}

\section{Convergence analysis of cyclic coordinatewise minimization for $f_2$} 
\label{sec:convergence:withlog}

\noindent
In this section, we consider the convergence behavior of the cyclic coordinatewise minimization 
algorithm applied to the function $f_2$ (Algorithm 2). It follows by assumption (A5)* and the 
convexity of $f_2$ that the set of optimal solutions of the minimization problem in 
(\ref{l1minpblm2}), denoted by $\mathcal{X}_\ell^*$, is non-empty. Since the negative logarithm 
function is convex, and $q$ is strictly convex, it follows by arguments very similar to those in 
\cite[Page 5]{Luo:Tseng:1989} that $\mathcal{X}_\ell^*$ is a convex set and that there exists 
${\bf t}^* \in \mathbb{R}^m$ such that $E {\bf x}^* = {\bf t}^*, \; \forall \; {\bf x}^* \in 
\mathcal{X}_\ell^*$. Let ${\bf d} ({\bf x}) = \nabla \{q(E {\bf x})\} = 2 E^T E {\bf x}$. We denote 
the $i^{th}$ entry of ${\bf d} ({\bf x})$ by $d_i ({\bf x})$. Let ${\bf d}^* := E^T E {\bf x}^* = 
E {\bf t}^*$. It follows that 
\begin{equation} \label{eqwlog1}
{\bf d} ({\bf x}^*) = {\bf d}^* \; \forall x^* \in \mathcal{X}_\ell^*. 
\end{equation}

\noindent
We now state two lemmas which will be important in understanding the coordinatewise 
minimization for the function $f_2$. 
\begin{lemma} \label{lemwlog1}
\begin{enumerate}[(a)]
\item Let $h(u) = au^2 + bu + c - \log u$ for $u > 0$. If $a > 0$, then $h(u)$ is uniquely minimized 
at $u^* = \frac{-b + \sqrt{b^2 + 8a}}{4a}$. 
\item Let $h(u) = au^2 + bu + c + \lambda |u|$ for $u \in \mathbb{R}$. If $a, \lambda > 0$, then 
$h(u)$ is uniquely minimized at $u^* = S_\lambda (-b)/2a$, where $S_\lambda$ is the 
soft-thresholding operator defined by $S_\lambda (x) = sign(x) (|x| - \lambda)_+$. 
\end{enumerate}
\end{lemma}

\noindent
{\it Proof} (a) Note that 
$$
\frac{d}{du} h(u) = 0 \Leftrightarrow 2au^2 + bu - 1 = 0. 
$$

\noindent
The result follows by noting that $u^*$ the only non-negative solution of the above 
equation, and that $h$ is a strictly convex function. 

\noindent
(b) The KKT conditions for the minimizing the strictly convex function $h$ are 
satisfied if and only if $u = 0$ if $|b| \leq \lambda$, and $2au + b + \lambda sign(u) = 0$ 
if $|b| > \lambda$, which in turn is satisfied if and only if $u = u^*$. \hfill$\Box$

\bigskip

\noindent
It is clear from Lemma \ref{lemwlog1} that the coordinatewise minimizers for $f_2$ (see 
(\ref{eqwlog3})) are uniquely defined and can be obtained in closed form. Note that the 
function $f_2 ({\bf x})$ takes the value infinity if $x_i = 0$ for $i$ belonging to a 
non-trivial subset of $S^c$. Hence, the KKT conditions for the convex minimization 
problem in (\ref{l1minpblm2}) imply that ${\bf x} \in \mathcal{X}_\ell^*$ if and only if 
\begin{eqnarray}
& & d_i ({\bf x})) = \frac{1}{x_i} \; \mbox{ for } i \in S^c, \label{kktwlog1}\\
& & d_i ({\bf x}) + \lambda sign(x_i) = 0 \; \mbox{ if } x_i \neq 0, i \in S, \label{kktwlog2}\\
& & |d_i ({\bf x})| \leq \lambda \; \mbox{ if } x_i = 0, i \in S. \label{kktwlog3} 
\end{eqnarray}

\noindent
The arguments in the proof of Lemma \ref{lem2} can be used to provide the following 
alternative characterization of the elements of $\mathcal{X}_\ell^*$. 
\begin{lemma} \label{lemwlog3}
${\bf x} \in \mathcal{X}_\ell^*$ if and only if 
\begin{eqnarray}
& & d_i ({\bf x})) = \frac{1}{x_i} \; \mbox{ for } i \in S^c, \label{kktwlog4}\\
& & x_i = sign(x_i - d_i ({\bf x})) \max(|x_i - d_i ({\bf x})| - \lambda, 0) \mbox{ for } i \in S. 
\label{kktwlog5} 
\end{eqnarray}
\end{lemma}

\noindent
Recall that $\{{\bf z}^r\}_{r \geq 0}$ is the sequence of iterates generated by 
Algorithm 2, and ${\bf x}^{r,i}$ is the appropriate coordinatewise minimizer defined in 
(\ref{eqwlog3}). It follows from arguments similar to those in the proof of Lemma \ref{lem2} 
that for $i \in S$, 
\begin{equation} \label{eqwlog4}
z^{r,i}_i  = sign(z^{r,i}_i - d_i ({\bf z}^{r,i})) \max(|z^{r,i}_i - d_i ({\bf z}^{r,i})| - \lambda, 0), 
\end{equation}

\noindent
and for $i \in S^c$ 
\begin{equation} \label{eqwlog5}
z^{r,i}_i = \frac{1}{d_i ({\bf z}^{r,i}))}. 
\end{equation}}

\noindent
As in Section \ref{sec:convergence}, we will establish a series of lemmas, which will ultimately lead 
us to the proof of Theorem \ref{thmwlog1}. Let 
$$
{\bf t}^{r,i} = E {\bf z}^{r,i} 
$$

\noindent
for all $r$ and all $0 \leq i \leq n$. By (\ref{eqwlog3}), it follows that 
\begin{equation} \label{eqwlog6}
f_2 ({\bf z}^{r,i}) \leq f_2 ({\bf z}^{r,i-1}) 
\end{equation}

\noindent 
for every $r$ and $1 \leq i \leq n$. It follows by assumption (A5)* that 
\begin{equation} \label{eqwlog6.1}
\{{\bf t}^{r,i}\}_{r \geq 0, 1 \leq i \leq n} \mbox{ is bounded}. 
\end{equation}

\noindent
By (\ref{eqwlog6}), the sequence $\{f_2 ({\bf z}^{r,i})\}_{r \geq 0}$ decreases to the same quantity, 
say $f^\infty$ for every $0 \leq i \leq n$. Since $\mathcal{X}_\ell^*$ is non-empty, it follows that 
$f^\infty > -\infty$. The next lemma shows that the sum of norm-square of the difference 
between successive iterates in $\{{\bf z}^r\}_{r \geq 0}$ is finite. Note that in 
Section \ref{sec:convergence}, we first needed to show that $\|{\bf z}^r - {\bf z}^{r+1}\|$ 
converges to zero (Lemma \ref{lem4}) to prove a similar result (Lemma \ref{lem5}). However, since 
we have to deal with the quadratic function $q(E{\bf x})$ as opposed to a general $g(E{\bf x})$ in 
this section, a direct argument is available. 
\begin{lemma} \label{lemwlog3.1}
$$
\sum_{r=0}^\infty \|{\bf z}^r - {\bf z}^{r+1}\|^2 < \infty. 
$$
\end{lemma}

\noindent
{\it Proof} For every $1 \leq i \leq n$, a second order Taylor series expansion along the $i^{th}$ 
coordinate leads to the following. 
\begin{equation} \label{eqwlog6.2}
g(E {\bf z}^{r,i-1}) - g(E {\bf z}^{r,i}) = d_i ({\bf z}^{r,i}) (z^{r,i-1}_i - z^{r,i}_i) + 
2 \|E_{\cdot i}\|^2 (z^{r,i-1}_i - z^{r,i}_i)^2. 
\end{equation}

\noindent
Fix $i \in S$ arbitrarily. Using exactly the same argument as in the proof of Lemma \ref{lem5} (b) for 
this case, we get that 
\begin{equation} \label{eqwlog7}
f({\bf z}^{r,i-1}) - f({\bf z}^{r,i}) \geq 2 \left( \min_{1 \leq j \leq n} \|E_{\cdot j}\|^2 \right) 
(z^{r,i-1}_i - z^{r,i}_i)^2. 
\end{equation}

\noindent
Fix $i \in S^c$ arbitrarily. By strict convexity of the negative logarithm function on $\mathbb{R}_+$, 
it follows that 
\begin{equation} \label{eqwlog8}
(-\log z^{r,i-1}_i) - (-\log z^{r,i}_i) \geq \left( -\frac{1}{z^{r,i}_i} \right) (z^{r,i-1}_i - z^{r,i}_i). 
\end{equation}

\noindent
It follows by (\ref{eqwlog5}), (\ref{eqwlog6.2}) and (\ref{eqwlog8}) that (\ref{eqwlog7}) is satisfied 
for every $i \in S^c$. Adding (\ref{eqwlog5}) over $i = 1,2, \cdots, n$, we obtain 
$$
 f({\bf z}^r) - f({\bf z}^{r+1}) \geq 2 \left( \min_{1 \leq j \leq n} \|E_{\cdot j}\|^2 \right) 
\sum_{i=1}^n (z^{r,i-1}_i - z^{r,i}_i)^2 = 2\left( \min_{1 \leq j \leq n} \|E_{\cdot j}\|^2 
\right) \|{\bf z}^r - {\bf z}^{r+1}\|^2. 
$$

\noindent
The result follows by noting that $f({\bf z}^r) \downarrow f^\infty > -\infty$ as $r \rightarrow 
\infty$ and that $\min_{1 \leq j \leq n} \|E_{\cdot j}\|^2 > 0$ as $E$ has no zero column. 
\hfill$\Box$ 

\medskip

\noindent
By Lemma \ref{lemwlog3.1}, it follows that $\|{\bf z}^r - {\bf z}^{r+1}\| \rightarrow 0$ as 
$r \rightarrow \infty$. We now establish a parallel version of Lemma \ref{lem5} for the 
problem at hand. 
\begin{lemma} \label{lemwlog4}
For every $0 \leq i \leq n$, 
\begin{equation} \label{eqwlog9}
\|{\bf t}^{r,i} - {\bf t}^*\| \rightarrow 0, 
\end{equation}

\noindent
as $r \rightarrow \infty$. 
\end{lemma}

\noindent
{\it Proof} By exactly the same set of arguments as in the proof of Lemma \ref{lem5}, there 
exists ${\bf t}^\infty \in \mathbb{R}^m$, and a subsequence $\mathcal{R}$ of $\mathbb{N}$ 
such that 
\begin{equation} \label{eqwlog10}
\{{\bf t}^{r,j}\}_{r \in \mathcal{R}} \rightarrow {\bf t}^\infty 
\end{equation}

\noindent
for every $0 \leq j \leq n$. Let $d^\infty = 2 E^T {\bf t}^\infty)$. It follows that 
\begin{equation} \label{eqwlog11}
\{d({\bf z}^{r,j}\}_{r \in \mathcal{R}} = \{2 E^T {\bf t}^{r,j}\}_{r \in \mathcal{R}} \rightarrow 
d^\infty 
\end{equation}

\noindent
as for every $0 \leq j \leq n$. Suppose $i \in S$. By repeating exactly the same arguments in the 
proof of Lemma \ref{lem5} in this case, we get the following. 
\begin{itemize}
\item If $|d^\infty_i| < \lambda$, then 
\begin{equation} \label{eqwlog12}
z^{r+1}_i = z^{r,i}_i = 0 
\end{equation}

\noindent
for large enough $r$. 
\item If $d^\infty_i = \lambda$, then 
\begin{equation} \label{eqwlog13}
z^{r+1}_i = z^{r,i}_i \leq 0 
\end{equation}

\noindent
for large enough $r$. 
\item If $d^\infty_i = -\lambda$, then 
\begin{equation} \label{eqwlog14}
z^{r+1}_i = z^{r,i}_i \geq 0 
\end{equation}

\noindent
for large enough $r$. 
\end{itemize}

\noindent
Since $\{f_2 ({\bf z}^{r,i})\}_{r \geq 0, 1 \leq i \leq n}$ is bounded above, it follows by 
assumption (A5)* that $\{{\bf z}^{r+1}\}_{r \in \mathcal{R}}$ is bounded (with the 
coordinates in $S^c$ uniformly bounded away from zero), and hence has at least one 
limit point. Let ${\bf z}^\infty$ denote any limit point of $\{{\bf z}^{r+1}\}_{r \in 
\mathcal{R}}$. It follows that 
\begin{equation} \label{eqwlog15}
E {\bf z}^\infty = {\bf t}^\infty \mbox{ and } d({\bf z}^\infty) = 2 E^T \nabla 
{\bf t}^\infty = d^\infty. 
\end{equation}

\noindent
It follows by (\ref{eqwlog5}), (\ref{eqwlog11}), (\ref{eqwlog12}), (\ref{eqwlog13}) and 
(\ref{eqwlog14}) that $z^\infty_j = 1/d^\infty_j$ if $j \in S^c$, $z^\infty_j = 0$ if $j 
\in S$ and $|d^\infty_j| < \lambda$, $z^\infty_j \leq 0$ if $j \in S$ and $d^\infty_j = 
\lambda$, $z^\infty_j \geq 0$ if $j \in S$ and $d^\infty_j = -\lambda$. It follows 
from Lemma \ref{lemwlog1} (b) that ${\bf z}^\infty \in \mathcal{X}_\ell^*$. It follows 
by Lemma \ref{lem1} that ${\bf t}^\infty = E {\bf z}^\infty = {\bf t}^*$. The result 
follows by noting that ${\bf t}^\infty$ is an arbitrarily chosen limit point of 
$\{{\bf t}^{r,i}\}_{r \geq 0}$. \hfill$\Box$

\noindent
It follows by Lemma \ref{lemwlog4} and the continuity of $g$ at ${\bf t}^*$ that 
\begin{equation} \label{eqwlog16}
{\bf d} ({\bf z}^{r,i}) \rightarrow {\bf d}^* 
\end{equation}

\noindent
as $r \rightarrow \infty$ for every $1 \leq i \leq n$. The next two lemmas show that the 
sequence of iterates $\{{\bf z}^r\}_{r \geq 0}$ approaches $\mathcal{X}_\ell^*$. 
\begin{lemma} \label{lemwlog5}
If $i \in S$, then 
$$
z^r_i - sign(z^r_i - d_i ({\bf z}^r)) \max(|z^r_i - d_i ({\bf z}^r)| - \lambda, 0) \rightarrow 0 
$$

\noindent
as $r \rightarrow \infty$. If $i \in S^c$, then 
$$
z^r_i \rightarrow \frac{1}{d_i^*}.  
$$

\noindent
as $r \rightarrow \infty$. 
\end{lemma}

\noindent
The proof of  the above lemma is provided in the appendix. As in Section \ref{sec:convergence}, for every ${\bf x} \in 
\mathcal{X}$, define the function $\phi$ as follows: 
$$
\phi ({\bf x}) = \min_{{\bf z}^* \in \mathcal{X}_\ell^*} \|{\bf x} - {\bf z}^*\|. 
$$

\noindent
Hence, $\phi({\bf x})$ is the distance of ${\bf x}$ from the closed convex set $\mathcal{X}_\ell^*$. 
\begin{lemma} \label{lemwlog6}
$$
\phi({\bf z}^r) \rightarrow 0 \mbox{ as } r \rightarrow \infty. 
$$
\end{lemma}

\noindent
{\it Proof} By (\ref{kktwlog1}), (\ref{kktwlog2}), (\ref{kktwlog3}) and the fact that $E {\bf z}^* 
= {\bf t}^*$ for every ${\bf z}^* \in \mathcal{X}_\ell^*$, it follows that $\mathcal{X}_\ell^*$ is the 
solution set of the linear system of equations given by 
\begin{eqnarray*}
& & E {\bf y} = {\bf t}^*, {\bf y} \in \mathcal{X},\\
& & y_i = 0 \mbox{ if } i \in S, |d_i^*| < \lambda,\\
& & y_i \leq 0 \mbox{ if } i \in S, d_i^* = \lambda,\\
& & y_i \geq 0 \mbox{ if } i \in S, d_i^* = -\lambda,\\
& & y_i = \frac{1}{d_i^*} \mbox{ if } i \in S^c. 
\end{eqnarray*}

\noindent
Note that the statements of Lemma \ref{lem8} apply exactly to the problem at hand for 
$i \in S$. Since $\mathcal{X}_\ell^*$ is non-empty, by Lemma \ref{lem3}, for sufficiently 
large $r$, there exists ${\bf y}^r \in \mathcal{X}_\ell^*$ such that 
\begin{eqnarray} 
\|{\bf z}^r - {\bf y}^r\| 
&\leq& \theta \left( \|E {\bf z}^r - {\bf t}^*\| + \sum_{i \in S, d_i^* = \lambda} (z^r_i)^+ + 
\sum_{i \in S, d_i^* = -\lambda} (-z^r_i)^+  + \sum_{i \in S, |d_i^*| < \lambda} |z^r_i| \right) + 
\nonumber\\
& & \theta \sum_{i \in S^c} \left| z^r_i - \frac{1}{d_i^*} \right| \nonumber\\
&=& \theta \left( \|E {\bf z}^r - {\bf t}^*\| + \sum_{i \in S^c} \left| z^r_i - \frac{1}{d_i^*} \right| 
\right), \label{eqwlog17}
\end{eqnarray}

\noindent
where $\theta$ is a constant only depending on $E$. The result follows by the definition of 
$\phi$, Lemma \ref{lemwlog4} and Lemma \ref{lemwlog5}. \hfill$\Box$ 

\medskip

\noindent
Let 
\begin{eqnarray*}
& & I_1^* := \{i \in S: d_i^* = \lambda\},\\
& & I_2^* := \{i \in S: d_i^* = -\lambda\},\\
& & I_3^* := \{i \in S: |d_i^*| < \lambda\}. 
\end{eqnarray*}

\noindent
By Lemma \ref{lem8} (recall that the statements of this lemma apply verbatim for $i \in S$) there 
exists an $r_0 > 0$ such that 
\begin{equation} \label{eqwlog18}
{\bf z}^r_{I_3^*} = {\bf 0}, \; ({\bf z}^r_{I_1^*})^+ = {\bf 0}, \; (-{\bf z}^r_{I_2^*})^+ = {\bf 0} 
\end{equation}

\noindent
for every $r \geq r_0$. Let $M = 2 E^T E$. By the assumption that $E$ has no zero column, 
it follows that $m_{ii} > 0$ for every $1 \leq i \leq n$. As in Section \ref{sec:convergence}, 
let $B$ denote the lower triangular portion of $M$, and $C = M - B$ denote the strictly upper 
triangular portion of $M$. Since $M$ is a positive semi-definite matrix with strictly positive 
diagonal entries, it follows that Lemma \ref{lem14} and Lemma \ref{lem15} hold with this 
choice of $M, B$ and $C$. 

\noindent
Let $I^* = I_1^* \cup I_2^*$, and 
$$
\beta = \max_{J \subseteq I^*} \sqrt{|J^c|} \left\{ \left( \frac{\tau_J \|(B_{JJ}^{-1}\| \|M_{JJ}\|}{1 - 
\rho_J} + \Delta + 1 \right) \|(B_{JJ})^{-1} B_{JJ^c}\| + \frac{\tau_J \|(B_{JJ})^{-1}\| \|M_{JJ}\|}{1 - 
\rho_J} \right\}. 
$$

\noindent
The next lemma is a parallel version of  analogous to Lemma \ref{lem16} for the problem at hand. 
\begin{lemma} \label{lemwlog7}
Consider any $J \subseteq I^*$. If for some two integers $s \geq t \geq r_0$ we have $z_i^r 
\neq 0$ for every $t+1 \leq r \leq s$ and $i \in J$, then, for any ${\bf z}^* \in \mathcal{X}_\ell^*$, 
there holds 
$$
\|{\bf z}_J^s - {\bf z}_J^*\| \leq \Delta \|{\bf z}_J^t - {\bf z}_J^*\| + \beta \max_{t \leq r \leq s} 
\|{\bf z}_{J^c}^r - {\bf z}_{J^c}^*\|_\infty. 
$$
\end{lemma}

\noindent
{\it Proof} Since $\Delta \geq 1$, it follows that the claim holds if $s = t$. Suppose $s > t \geq 
r_0$. Fix any $r \in \{t, \cdots, s-1\}$ and $i \in I^*$. Recall that $q({\bf y}) = {\bf y}^T {\bf y}$. 
By using exactly the same arguments as in the beginning of the proof of Lemma \ref{lem16}, it 
follows that 
\begin{eqnarray*}
0 = d_i ({\bf z}^{r,i}) - d_i^* = 2 E_{\cdot i}^T E {\bf z}^{r,i} - 2 E_{\cdot i}^T {\bf t}^* = 
E_{\cdot i}^T \nabla^2 q({\bf t}^*)   (E {\bf z}^{r,i} - {\bf t}^*). 
\end{eqnarray*}

\noindent
The result now follows by using exactly the same argument as in the proof of 
\cite[Lemma 9]{Luo:Tseng:1989} (starting from \cite[Page 12, Line -5]{Luo:Tseng:1989} to 
the end of the proof, replacing $E$ by $E^T$, and $w_J^r$ by ${\bf 0}$ throughout). 
\hfill$\Box$ 

\medskip

\noindent
Let $\sigma_0 := 1$ and 
$$
\sigma_k = \Delta + 3+ \beta + (\beta + 1) \sigma_{k-1} \; k = 1,2, \cdots, n. 
$$

\noindent
It follows from the above definition that $\sigma_k \geq 1$ for every $1 \leq k \leq n$, and is 
monotonically increasing with $k$. 

Fix $\delta > 0$ arbitrarily. Note that by (\ref{kktwlog1}), $z_i^* = 1/d_i^*$ for every $i \in 
S^c$ and every ${\bf z}^* \in \mathcal{X}_\ell^*$. By Lemma \ref{lemwlog3.1} and 
Lemma \ref{lemwlog6} there exists $r_1 > 0$ such that for every $r \geq r_1$, 
\begin{eqnarray}
\phi ({\bf z}^r) \leq \delta, \label{eqwlog19}\\
\|{\bf z}^{r+1} - {\bf z}^r\| \leq \delta, \label{eqwlog20}\\
\|{\bf z}_{S^c}^r - {\bf z}_{S^c}^*\| \leq \delta, \mbox{ for every } {\bf z}^* \in 
\mathcal{X}_\ell^*. \label{eqwlog21}
\end{eqnarray}

\noindent
The next three lemmas are parallel versions of Lemma \ref{lem17}, Lemma \ref{lem18} and 
Lemma \ref{lem19} respectively. The proofs of these lemmas follow by repeating the proofs 
of Lemma \ref{lem17}, Lemma \ref{lem18} and Lemma \ref{lem19} verbatim, with the 
following exceptions: replace $n$ by $|S|$ throughout, replace $i \notin J$ by $i \notin 
S \setminus J$, and replace $\mu$ by $0$. 
\begin{lemma} \label{lemwlog8}
Fix $k \in \{1,2, \cdots, |S|\}$ arbitrarily. If for some nonempty $J \subset I^*$, and some 
intergers $t' > t \geq \max (r_0, r_1)$, we have 
\begin{eqnarray}
& & |z_i^t| > \sigma_k \delta, \; \forall i \in J, \label{eqwlog22}\\
& & |z_i^r| \leq \sigma_{k-1} \delta, \; \forall i \notin J, \setminus S \forall r = t, t+1, \cdots, t' - 1 
\label{eqwlog23}, 
\end{eqnarray}

\noindent
then the following hold: 
\begin{enumerate}[(a)]
\item $|z_i^{t'}| > \sigma_{k-1} \delta$ for every $i \in J$. 
\item There exists an ${\bf z}^* \in \mathcal{X}_\ell^*$ such that 
$$
\|{\bf z}^r - {\bf z}^*\|_\infty \leq \sigma_k \delta, \; \forall r = t, t+1, \cdots, t'-1. 
$$
\end{enumerate}
\end{lemma}

\begin{lemma} \label{lemwlog9}
Fix $k \in \{1,2, \cdots, |S|\}$ arbitrarily. If for some $J \subseteq I^*$ with $|J| \geq |I^*| - 
k + 1$ and some interger $t > \max(r_0, r_1)$ we have 
\begin{eqnarray}
& & |z_i^t| > \sigma_k \delta, \; \forall i \in J, \label{eqwlog24}\\
& & |z_i^t| \leq \sigma_{k-1} \delta, \; \forall i \notin S \setminus J, \label{eqwlog25} 
\end{eqnarray}

\noindent
then there exists an ${\bf z}^* \in \mathcal{X}_\ell^*$ and a $\bar{t} \geq t$ satisfying 
\begin{equation} \label{eqwlog26}
\|{\bf z}^r - {\bf z}^*\|_\infty \leq \sigma_k \delta, 
\end{equation}

\noindent
for every $r \geq \bar{t}$. 
\end{lemma}

\begin{lemma} \label{lemwlog10}
For any $\delta > 0$, there exists an ${\bf z}^* \in \mathcal{X}_\ell^*$ and $\hat{r} > 0$ such 
that 
\begin{equation} \label{eqwlog27}
\|{\bf z}^r - {\bf z}^*\|_\infty \leq \sigma_{|S|} \delta + \delta, 
\end{equation}

\noindent
for every $r \geq \hat{r}$. 
\end{lemma}

\noindent
We can now prove Theorem \ref{thmwlog1} by repeating the arguments at the end of 
Section \ref{sec:convergence} (after the proof of Lemma \ref{lem19}) verbatim.

\section{Applications} \label{sec:applications}

\noindent
In this section, we demonstrate the utility of Theorem \ref{thm1} and 
Theorem \ref{thmwlog1}. In particular, we use these results to establish convergence of two 
commonly used cyclic coordinatewise descent algorithms: one arising in high dimensional 
covariance estimation in the context of graphical models, and another arising in high dimensional 
logistic regression. 

\subsection{Convergence of a pseudo likelihood based algorithm for graphical model selection}

\noindent
The CONCORD algorithm, introduced in Khare et al. \cite{KOR:2014:CONCORD}, is a 
sparse inverse covariance estimation algorithm, which uses cyclic coordinatewise 
minimization to minimize the function 
\begin{equation} \label{eq54}
Q_{con} (\Omega) = \sum_{i=1}^p - \log \omega_{ii} + \frac{1}{2} \sum_{i=1}^p 
\Omega_{\cdot i}^T \hat{\Sigma} \Omega_{\cdot i} + \lambda \sum_{1 \leq i < j \leq p} 
|\omega_{ij}|, 
\end{equation}

\noindent
subject to the constraint that $\Omega = ((\omega_{ij}))_{1 \leq i,j \leq p}$ is a $p 
\times p$ symmetric matrix with non-negative diagonal entries. Here $\Omega_{\cdot i}$ 
denotes the $i^{th}$ column of $\Omega$, $p$ is a fixed positive integer, $\lambda 
> 0$ is a fixed positive real number, and $\hat{\Sigma}$ is the (observed) sample 
covariance matrix of $n$ i.i.d. observations from a $p$-variate distribution. Hence 
$\hat{\Sigma}$ is positive semi-definite. Let $\Sigma = \Omega^{-1}$ denote the (unknown) true covariance matrix for the 
underlying $p$-variate distribution. The CONCORD algorithm provides a sparse estimate of the inverse covariance matrix 
$\Omega$ by minimizing the objective function $Q_{con}$. Models which induce sparsity in the inverse covariance matrix 
are known as concentration graphical models, and have gained popularity in statistics, machine learning etc.  

As with any sparse covariance estimation algorithm, the CONCORD algorithm is particularly developed to tackle 
high-dimensional settings, i.e., settings where $p$ is much larger than $n$.
The function $Q_{con} (\Omega)$ is a convex function of $\Omega$, 
but is not necessarily strictly convex if $n < p$, as 
the matrix $\hat{\Sigma}$ is singular in this case. 

Other pseudo likelihood based sparse inverse covariance estimation algorithms in the literature (see 
\cite{KOR:2014:CONCORD} for a list of references) also provide sparse estimates of $\Omega$ via cyclic 
coordinatewise minimization for objective functions which are different from 
$Q_{con}$. However, there are no convergence guarantees for the corresponding 
algorithms. In fact, as shown in \cite{KOR:2014:CONCORD}, it is easy to find (non-pathological)
examples where some of these algorithms do not converge. On the other hand, as shown below, the results in this paper can be used to establish convergence of the CONCORD algorithm. 

Note that the output produced by the CONCORD algorithm is not guaranteed to be positive 
definite (same is true for the algorithms in \cite{Peng:2009, Yu:2008}). However, the focus here 
is to estimate the sparsity pattern in $\Omega$, i.e., model selection. If needed, a positive definite 
version with the estimated sparsity pattern can be constructed using standard approaches (see 
Khare et al. \cite{KOR:2014:CONCORD} for a discussion). 

We first provide a lemma which will be useful in our convergence proof. 
\begin{lemma} \label{lem20}
Let $A$ be a $k \times k$ positive semi-definite matrix with $A_{kk} > 0$, and $\lambda$ be a 
positive constant. Consider the function 
$$
h({\bf x}) = -\log x_k + {\bf x}^T A {\bf x} + \lambda \sum_{i=1}^{k-1} |x_j| 
$$

\noindent
defined on $\mathbb{R}^{k-1} \times \mathbb{R}_+$. Then, there exist positive constants $a_1$
and $a_2$ (depending only on $\lambda$ and $A$), such that 
$$
h({\bf x}) \geq a_1 x_k - a_2 
$$

\noindent
for every ${\bf x} \in \mathbb{R}^{k-1} \times \mathbb{R}_+$. 
\end{lemma}

\noindent
{\it Proof} Let ${\bf x}_{-k} := (x_i)_{1 \leq i \leq k-1}$, and 
$$
A = \left[ \begin{matrix}
A_1 & {\bf b} \cr
{\bf b}^T & A_{kk} 
\end{matrix}
\right]. 
$$

\noindent
Since $A$ is positive semi-definite, and $A_{kk} > 0$, it follows that 
\begin{equation} \label{eq54.1}
{\bf x}^T A {\bf x} = A_{kk} \left( x_k + \frac{{\bf b}^T {\bf x}_{-k}}{A_{kk}} \right)^2 + 
{\bf x}_{-k}^T \left( A_1 - \frac{1}{A_{kk}} {\bf b} {\bf b}^T \right) {\bf x}_{-k} \geq A_{kk} 
\left( x_k + \frac{{\bf b}^T {\bf x}_{-k}}{A_{kk}} \right)^2. 
\end{equation}

\noindent
Note that for any $c > 0$, the function $cy - \log y$ is minimized at $y = \frac{1}{c}$. Hence, 
for every $c > 0$ and $y > 0$, we get that $cy - \log y \geq 1 + \log c$. If ${\bf b} = 0$, then 
it follows by (\ref{eq54.1}) and the definition of $h({\bf x})$ that 
\begin{equation} \label{eq54.2}
h({\bf x}) \geq - \log x_k + A_{kk} x_k^2 \geq - \log x_k + 2 A_{kk} x_k - A_{kk} \geq A_{kk} x_k 
+ 1 + \log A_{kk} - A_{kk}. 
\end{equation}
 
\noindent
Hence the result holds if ${\bf b} = {\bf 0}$. 

If ${\bf b} \neq {\bf 0}$, then $\|{\bf b}\|_\infty > 0$. Since $|{\bf b}^T {\bf x}_{-k}| \leq 
\|{\bf b}\|_\infty \sum_{i=1}^{k-1} |x_i|$, it follows by (\ref{eq54.2}) and the definition of 
$h({\bf x})$ that 
\begin{equation} \label{eq54.3}
h({\bf x}) \geq - \log x_k + A_{kk} \left( x_k + \frac{{\bf b}^T {\bf x}_{-k}}{A_{kk}} \right)^2 + 
\frac{\lambda A_{kk}}{\|{\bf b}\|_\infty} \left| \frac{{\bf b}^T {\bf x}_{-k}}{A_{kk}} \right|. 
\end{equation}

\noindent
It follows by Lemma \ref{lemwlog1} (b) that for every $x_k > 0$ and $\widetilde{\lambda} > 
0$, the function $h(y) = (x_k + y)^2 + \widetilde{\lambda} |y|$ is minimized at $y = -\left( 
x_k - \frac{\widetilde{\lambda}}{2} \right)_+$. It follows from (\ref{eq54.3}) that 
\begin{eqnarray}
h({\bf x}) 
&\geq& - \log x_k + A_{kk} \left( x_k -\left( x_k - \frac{\lambda}{2 \|{\bf b}\|_\infty} 
\right)_+ \right)^2 + \frac{\lambda A_{kk}}{\|{\bf b}\|_\infty} \left( x_k - \frac{\lambda}{2 
\|{\bf b}\|_\infty} \right)_+ \nonumber\\
&\geq& -\log x_k + \min \left( A_{kk} x_k^2, \frac{\lambda A_{kk}}{\|{\bf b}\|_\infty} \left( 
x_k - \frac{\lambda}{2 \|{\bf b}\|_\infty} \right) \right) \label{eq54.4}
\end{eqnarray}

\noindent
The result follows from (\ref{eq54.4}), the fact that $x_k^2 \geq 2 x_k - 1$, and the fact that 
$cy - \log y \geq 1 + \log c$ (for $c = 1$ and $c = \lambda A_{kk}/(2 \|{\bf b}\|_\infty))$. 
\hfill$\Box$

\medskip

\noindent
The following theorem establishes the convergence of the CONCORD algorithm by using 
Theorem \ref{thmwlog1}. 
\begin{thm} \label{thm2}
If the diagonal entries of $\hat{\Sigma}$ are strictly positive, then the sequence of iterates generated by the cyclic 
coordinatewise minimization algorithm for $Q_{con}$ converges. 
\end{thm}

\noindent
{\it Proof} We will show that the minimization problem in (\ref{eq54}) is a special case 
of the minimization problem in (\ref{l1minpblm2}), and satisfies assumption (A5)*. Applying 
Theorem \ref{thmwlog1} yields the proof of convergence of the CONCORD algorithm. 

Let ${\bf y} = {\bf y} (\Omega) \in \mathbb{R}^{p^2}$ denote a vectorized version of 
$\Omega$ obtained by shifting the corresponding diagonal entry at the bottom of each 
column of $\Omega$, and then stacking the columns on top of each other. More precisely, if 
$P^i$ is the $p \times p$ permutation matrix such that $P^i {\bf z} = (z_1, \cdots, z_{i-1}, 
z_{i+1}, \cdots, z_p, z_i)$ for every ${\bf z} \in \mathbb{R}^p$, then 
$$
{\bf y} = {\bf y} (\Omega) = ((P^1 \Omega_{\cdot 1})^T, (P^2 \Omega_{\cdot 2})^T, \cdots, 
(P^p \Omega_{\cdot p})^T )^T. 
$$

\noindent
Note that since $\Omega$ is symmetric, $\omega_{ij} = \omega_{ji}$ for every $1 \leq i < j  
\leq p$. Let ${\bf x} = {\bf x} ({\Omega}) \in \mathbb{R}^{\frac{p(p+1)}{2}}$ be the 
symmetric version of ${\bf y}$, obtained by removing all $\omega_{ij}$ with $i > j$ 
from ${\bf y}$. More precisely, 
$$
{\bf x} = {\bf x} (\Omega) = (\omega_{11}, \omega_{12}, \omega_{22}, \cdots, 
\omega_{1p}, \omega_{2p}, \cdots, \omega_{pp})^T. 
$$

\noindent
Let $\tilde{P}$ be the $p^2 \times \frac{p(p+1)}{2}$ matrix such that every entry of $\tilde{P}$ 
is either $0$ or $1$, exactly one entry in each row of $\tilde{P}$ is equal to $1$, and ${\bf y} 
= \tilde{P} {\bf x}$. Let $\tilde{S}$ be a $p^2 \times p^2$ block diagonal matrix with $p$ 
diagonal blocks, and the $i^{th}$ diagonal block is equal to $\tilde{S}^i := \frac{1}{2} 
P^i \hat{\Sigma} (P^i)^T$. It follows that 
\begin{eqnarray}
\frac{1}{2} \sum_{i=1}^p \Omega_{\cdot i}^T \hat{\Sigma} \Omega_{\cdot i} = \frac{1}{2} 
\sum_{i=1}^p \Omega_{\cdot i}^T (P^i)^T P^i \hat{\Sigma} (P^i)^T P^i \Omega_{\cdot i} 
&=& \frac{1}{2} \sum_{i=1}^p (P^i \Omega_{\cdot i})^T (P^i \hat{\Sigma} (P^i)^T) (P^i 
\Omega_{\cdot i}) \nonumber\\
&=& {\bf y}^T \tilde{S} {\bf y} \nonumber\\
&=& {\bf x}^T \tilde{P}^T \tilde{S} \tilde{P} {\bf x}. \label{eq55} 
\end{eqnarray}

\noindent
Note that for every $1 \leq i \leq p$, the matrix $\tilde{S}^i = \frac{1}{2} P^i \hat{\Sigma} 
(P^i)^T$ is positive semi-definite. Let $\tilde{S}^{1/2}$ denote the $p^2 \times p^2$ block 
diagonal matrix with $p$ diagonal blocks, such that the $i^{th}$ diagonal block is given by 
$(\tilde{S}^i)^{1/2}$. Let $E = \tilde{S}^{1/2} \tilde{P}$. It follows by (\ref{eq55}) that 
\begin{equation} \label{eq56}
\frac{1}{2} \sum_{i=1}^p \Omega_{\cdot i}^T \hat{\Sigma} \Omega_{\cdot i} = (E {\bf x})^T 
(E {\bf x}). 
\end{equation}

\noindent
By the definition of ${\bf x} (\Omega)$, we obtain 
\begin{equation} \label{eq59}
\omega_{ii} = x_{\frac{i(i+1)}{2}} 
\end{equation} 

\noindent
for every $1 \leq i \leq p$. Let 
$$
S = \left\{ j: \; 1 \leq j \leq \frac{p(p+1)}{2}, \; j \neq \frac{i(i+1)}{2} \mbox{ for any } 1 
\leq i \leq p \right\}, 
$$ 

\noindent
and 
$$
\mathcal{X} = \{{\bf x} \in \mathbb{R}^{\frac{p(p+1)}{2}}: x_j \geq 0 \mbox{ for every } j 
\in S^c\}. 
$$

\noindent
It follows by (\ref{eq54}), (\ref{eq56}) and (\ref{eq59}) that the CONCORD algorithm can be viewed 
as a cyclic coordinatewise minimization algorithm to minimize the function 
\begin{equation} \label{eq60}
Q_{con} ({\bf x}) =  {\bf x}^T E^T E {\bf x} - \sum_{i \in S^c} \log x_i +  \lambda \sum_{i \in S} 
|x_i|, 
\end{equation}

\noindent
subject to ${\bf x} \in \mathcal{X}$. For any $1 \leq i \leq p(p+1)/2$, there exist $1 \leq k,l 
\leq p$ such that $x_i = \omega_{kl}$. Note that $\|E_{\cdot i}\|^2 = \frac{\hat{\Sigma}_{kk} + 
\hat{\Sigma}_{ll}}{2} > 0$. In order to verify assumption (A5)*, we consider the set the set 
$R_\xi = \{{\bf x}: Q_{con} ({\bf x}) \leq \xi\}$. Recall by (\ref{eq54}) 
that 
\begin{equation} \label{eq60.1}
Q_{con} ({\bf x}) = Q_{con} ({\bf x} (\Omega)) = \sum_{i=1}^p \left\{ - \log \omega_{ii} + 
\frac{1}{2} \Omega_{\cdot i}^T \hat{\Sigma} \Omega_{\cdot i} + \frac{\lambda}{2} \sum_{1 
\leq j \neq i \leq p} |\omega_{ij}| \right\}. 
\end{equation}

\noindent
It follows by applying Lemma \ref{lem20} for every $1 \leq i \leq p$ in (\ref{eq60.1}) that 
there exist positive constants $a_1$ and $a_2$ (depending only on $\hat{\Sigma}$ and 
$\lambda$) such that 
\begin{equation} \label{eq60.2}
Q_{con} ({\bf x}) \geq a_1 \sum_{i=1}^p \omega_{ii} - a_2 = a_1 \sum_{i \in S^c} x_i - a_2. 
\end{equation}

\noindent
Hence, if ${\bf x} \in R_\xi$, then 
\begin{equation} \label{eq60.3}
x_i \leq (\xi + a_2)/a_1 \leq \widetilde{\xi} 
\end{equation}

\noindent
for every $i \in S^c$, where $\widetilde{\xi} = (|\xi| + a_2)/a_1$. It also follows by the 
definition of $Q_{con}$ that if ${\bf x} \in R_\xi$, then $-\sum_{i \in S^c} \log x_i < \xi$. 
Hence $\prod_{i \in S^c} x_i > e^{-\xi}$. It follows by (\ref{eq60.3}) that 
\begin{equation} \label{eq60.4}
x_i > \frac{e^{-\xi}}{\widetilde{\xi}^{p-1}} 
\end{equation}

\noindent
for every $i \in S^c$. It follows by (\ref{eq60.3}) that if ${\bf x} \in R_\xi$, then 
\begin{equation} \label{eq60.5}
\sum_{i \in S} |x_i| \leq \xi + \sum_{i \in S^c} \log x_i \leq \xi + p \log \widetilde{\xi}. 
\end{equation}

\noindent
It follows by (\ref{eq60.3}), (\ref{eq60.4}) and (\ref{eq60.5}) that $Q_{con} ({\bf x})$ satisfies 
assumption (A5)*. Combining this with the continuity and convexity of $Q_{con}$, it 
follows that the set 
$$
\mathcal{X}_\ell^* = \{{\bf x} \in \mathcal{X}: \; Q({\bf x}) < \infty, \; Q({\bf x}^*) \leq 
Q({\bf x}) \mbox{ for every } {\bf x} \in \mathcal{X}\}. 
$$

\noindent
is non-empty. Hence, assumption (A5) holds. It follows by Theorem \ref{thmwlog1} and Remark 
\ref{rem1} that the sequence of iterates produced by the CONCORD algorithm converges. 
\hfill$\Box$

\medskip

\noindent
{\it Remark} Note that if $n \geq 2$, and none of the underlying marginal distributions is degenerate, 
then the diagonal entries of $\hat{\Sigma}$ are strictly positive, and the assumption in 
Theorem \ref{thm2} is immediately satisfied. 

\subsection{Convergence of $\ell_1$ minimization for logistic regression}

\noindent
Let $Y_1, Y_2, \cdots, Y_N$ denote independent random variables taking values in $\{-1,1\}$, and 
$\{{\bf z}^i\}_{i=1}^N$ be a collection of vectors in $\mathbb{R}^p$ such that 
$$
P(Y_i = y \mid {\bf z}^i) = \frac{1}{1 + e^{-y_i {\boldsymbol \beta}^T {\bf z}^i}} 
$$

\noindent
for every $1 \leq i \leq N$. The above statistical model is known as the logistic regression model, 
and the objective is to estimate the parameter ${\boldsymbol \beta} \in \mathbb{R}^p$. However, 
in many modern applications, the number of observations $N$ is much less than the number of 
parameters $p$. To tackle such a situation, Shevade and Keerthi \cite{Shevade:Keerthi:2003} 
(see also \cite{FHT:2010, KKB:2007, LLAN:2006, Park:Hastie:2007, Yun:Toh:2011}) propose 
estimating ${\boldsymbol \beta}$ by minimizing the following objective function: 
\begin{equation} \label{eq61}
Q_{logit} ({\boldsymbol \beta}) = \sum_{i=1}^N \log \left( 1 + e^{-y_i {\boldsymbol \beta}^T 
{\bf z}^i} \right) + \lambda \sum_{j=1}^p |\beta_j|. 
\end{equation}

\noindent
Here $y_1, y_2, \cdots, y_N$ denote the observed values of $Y_1, Y_2, \cdots, Y_N$ respectively, 
and $\lambda > 0$ is fixed. The purpose of adding the $\ell_1$ penalty term $\lambda 
\sum_{j=1}^p |\beta_j|$ is to induce sparsity in the parameter estimate. Consider the function 
$g: \mathbb{R}^N \rightarrow \mathbb{R}$ defined by 
\begin{equation} \label{eq61.1}
g({\boldsymbol \eta}) = \sum_{i=1}^N \log \left( 1 + e^{-y_i \eta_i} \right). 
\end{equation}

\noindent
Since 
$$
\frac{\partial^2}{\partial \eta_i^2} \log \left( 1 + e^{-y_i \eta_i} \right) = \frac{y_i^2 e^{-y_i 
\eta_i}}{(1 + e^{-y_i \eta_i})^2} > 0 
$$

\noindent
for every $1 \leq i \leq N$, it follows that $g$ is a strictly convex function. Let $X$ denote the 
$N \times p$ matrix with $i^{th}$ row given by $({\bf z}^i)^T$. It follows by (\ref{eq61}) that 
\begin{equation} \label{eq62}
Q_{logit} (\boldsymbol \beta) = g(X {\boldsymbol \beta}) + \lambda \sum_{j=1}^p |\beta_j|. 
\end{equation}

\noindent
Note that in a typical high-dimensional setting, we have $N < p$. Hence, the matrix $X$ is 
singular, and consequently the function $Q_{logit}$ is {\it not necessarily strictly convex}. 

Shevade and Keerthi \cite[Page 2248]{Shevade:Keerthi:2003} propose using cyclic coordinatewise 
minimization for minimizing $Q_{logit}$. We note that the final algorithm that they present (see 
\cite[Page 2249]{Shevade:Keerthi:2003}) is a variant where at each iteration, the ``best 
coordinate" is chosen according to an appropriate criterion, and the function is minimized with 
respect to the chosen coordinate (keeping all the other coordinates fixed). Note that the 
minimizer with respect to a single coordinate cannot be obtained in closed form in this situation. 
However, such a minimization involves a convex function on a susbet of $\mathbb{R}$, and 
numerical methods can be used to obtain the required minimizer accurately in a few steps. 
In particular, the authors in \cite{Shevade:Keerthi:2003} use a combination of Newton-Raphson 
and bisection methods. To conclude, coordinatewise minimization (cyclic or the variant 
approach described above) is a viable approach for this problem, and has been used in 
applications. 

It is claimed in \cite{Shevade:Keerthi:2003} that convergence follows from 
\cite[Prop. 4.1, Chap. 3]{Bertsekas:Tsitsiklis:1989}. However, the result 
\cite[Prop. 4.1, Chap. 3]{Bertsekas:Tsitsiklis:1989} states that if $F({\bf x})$ is a 
convex function, and the sequence $\{{\bf x}^t\}_{t \geq 0}$ is generated 
by using 
$$
{\bf x}^{t+1} = \mbox{arg} \min \left\{ F({\bf x}) + \frac{1}{2c_t} \|{\bf x} - {\bf x}^t\|_2^2 
\right\}, 
$$

\noindent
where $\lim \inf_{t \rightarrow \infty} c_t > 0$; then $\{{\bf x}^t\}_{t \geq 0}$ converges to 
a global minimizer of $F({\bf x})$. Hence, to the best of our understanding, this result in 
not applicable to coordinatewise minimization in the current setting. 

We now show that Theorem \ref{thm1} can be used to provide a proof of convergence of the 
cyclic coordinatewise minimization algorithm for minimizing $Q_{logit}$. 
\begin{thm} \label{thm3}
If the matrix $X$ has no zero columns, then the sequence of iterates generated by the cyclic 
coordinatewise minimization algorithm for $Q_{logit}$ converges. 
\end{thm}

\noindent
{\it Proof} Consider the function $g$ defined in (\ref{eq61.1}). Note that $g$ is non-negative 
and $C_g = \mathbb{R}^N$. It follows by (\ref{eq62}) that the minimization problem for 
$Q_{logit}$ is a special case of the minimization problem in (\ref{l1minpblm1}) with $m = N$, 
$n = p$ and $E = X$, and that assumptions (A1)-(A4) are satisfied. 
Also, if $Q_{logit} ({\boldsymbol \beta}) \leq \xi$, it follows that $|\beta_j| \leq 
\xi/\lambda$ for every $1 \leq j \leq p$. Hence, the set $\{{\boldsymbol \beta}: 
Q_{logit} ({\boldsymbol \beta}) \leq \xi\}$ is a bounded set for every $\xi \in \mathbb{R}$. It 
follows that $Q({\boldsymbol \beta})$ satisfies assumption (A5). The result now follows by 
Theorem \ref{thm1}. \hfill$\Box$

\section*{Appendix}

\noindent
{\it Proof of Lemma \ref{lem2}} First, let us assume that $x \in \mathcal{X}^*$. Then (\ref{kkt1}), (\ref{kkt2}) and 
(\ref{kkt3}) hold. Hence (\ref{kkt4}) holds automatically. Suppose $i \in S$ and $x_i \neq 0$. By 
(\ref{kkt2}), it follows that $d_i ({\bf x}) + \lambda sign(x_i) = 0$. Hence, 
$$
|x_i - d_i ({\bf x})| = |x_i + \lambda sign(x_i)| = |x_i| + \lambda. 
$$

\noindent
Since $\lambda, |x_i| > 0$, we obtain $sign(x_i + \lambda sign(x_i)) = sign(x_i)$. It follows that 
$$
sign(x_i - d_i ({\bf x})) \max(|x_i - d_i ({\bf x})| - \lambda, 0) = sign(x_i) \max(|x_i| + \lambda - 
\lambda, 0) = x_i. 
$$

\noindent
Suppose $i \in S$ and $x_i = 0$. By (\ref{kkt3}), it follows that 
$$
sign(x_i - d_i ({\bf x})) \max(|x_i - d_i ({\bf x})| - \lambda, 0) = sign(-d_i ({\bf x})) \max(|d_i ({\bf x})| - 
\lambda, 0) = 0. 
$$

\noindent
This establishes (\ref{kkt5}). 

Now, let us assume that (\ref{kkt4}) and (\ref{kkt5}) hold. Hence (\ref{kkt1}) holds automatically. 
Suppose $i \in S$ and $x_i \neq 0$. We consider three cases. 
\begin{enumerate}
\item If $x_i  = d_i ({\bf x})$, then $x_i = 0$ by (\ref{kkt5}), which is a contradiction. 
\item If $x_i > d_i ({\bf x}$, then $x_i = \max(x_i - d_i ({\bf x}) - \lambda, 0)$ by (\ref{kkt5}). 
Since $x_i \neq 0$, it follows that $x_i > 0$ and $d_i ({\bf x}) = - \lambda = - \lambda sign(x_i)$. 
\item If $x_i < d_i ({\bf x})$, then $x_i = min(x_i + \lambda - d_i ({\bf x}), 0)$ by (\ref{kkt5}). 
Since $x_i \neq 0$, it follows that $x_i < 0$ and $d_i ({\bf x}) = \lambda = - \lambda sign(x_i)$. 
\end{enumerate}

\noindent
Hence (\ref{kkt2}) holds. Suppose $i \in S$ and $x_i = 0$. Then $\max(|d_i ({\bf x})| - \lambda, 0) 
= 0$ by (\ref{kkt5}). It follows that $|d_i ({\bf x})| \leq \lambda$. Hence, (\ref{kkt3}) holds. It follows 
that $x \in \mathcal{X}^*$. \hfill$\Box$ 

\bigskip

\noindent
{\it Proof of Lemma \ref{lem16}} Since $\Delta \geq 1$, it follows that the claim holds if $s = t$. Suppose $s > t \geq 
r_0$. Fix any $r \in \{t, \cdots, s-1\}$ and $i \in I^*$. Note that $x_i^{r+1} \neq 0$. If $i \in 
I_1^*$, it follows from (\ref{eq28}) that $x_i^{r+1} < 0$. Hence, by (\ref{eq4}), we obtain 
$d_i ({\bf x}^{r,i}) = \lambda = d_i^*$. If $i \in I_2^*$, it follows from (\ref{eq28}) that 
$x_i^{r+1} > 0$. Hence, by (\ref{eq4}), we obtain $d_i ({\bf x}^{r,i}) = -\lambda = d_i^*$. If $i 
\in I_4^*$, it follows from (\ref{eq28}) that $x_i^{r+1} > 0$. Hence, by (\ref{eq5}), we obtain 
$d_i ({\bf x}^{r,i}) = 0 = d_i^*$. In either case, it follows that 
\begin{eqnarray*}
0 
&=& d_i ({\bf x}^{r,i}) - d_i^*\\
&=& d_i ({\bf x}^{r,i}) - d_i ({\bf x}^*)\\
&=& E_{\cdot i}^T (\nabla g(E {\bf x}^{r,i}) - g(E {\bf x}^*))\\
&=& E_{\cdot i}^T \nabla^2 g(E {\bf x}^*) (E {\bf x}^{r,i} - E {\bf x}^*) + O \left( \|E {\bf x}^{r,i} 
- E {\bf x}^*\|^2 \right). 
\end{eqnarray*}

\noindent
The result now follows by using exactly the same argument as in the proof of 
\cite[Lemma 9]{Luo:Tseng:1989} (starting from \cite[Page 12, Line -5]{Luo:Tseng:1989} to 
the end of the proof, after replacing $E$ by $E^T$ throughout). \hfill$\Box$ 

\bigskip

\noindent
{\it Proof of Lemma \ref{lem17}} Let ${\bf x}^*$ be any element of $\mathcal{X}^*$ satisfying $\phi ({\bf x}^t) = 
\|{\bf x}^t - {\bf x}^*\|$. Hence, 
\begin{equation} \label{eq44}
\|{\bf x}^t - {\bf x}^*\| \leq \delta. 
\end{equation}

\noindent
By (\ref{eq43}), if $i \notin J$, then 
$$
|x_i^*| \leq |x_i^t| + \|{\bf x}^t - {\bf x}^*\| \leq \sigma_{k-1} \delta + \delta. 
$$

\noindent
It follows by (\ref{eq43}) that 
$$
||x_i^r| - |x_i^*|| \leq \sigma_{k-1} \delta + \delta, 
$$

\noindent
for every $t \leq r \leq t'-1$. Since $t \geq r_0$ it follows that $x_i^r$ and $x_i^*$ are either 
both non-positive or both non-negative for every $r \geq t$. Hence $||x_i^r| - |x_i^*|| = 
|x_i^r - x_i^*|$ for every $r \geq t$, which implies that 
\begin{equation} \label{eq45}
|x_i^r - x_i^*| \leq \sigma_{k-1} \delta + \delta, 
\end{equation}

\noindent
for every $t \leq r \leq t'-1$. We now claim that 
\begin{equation} \label{eq46}
|x_i^r| > \sigma_{k-1} \delta + \delta, 
\end{equation}

\noindent
for every $i \in J$ and $t \leq r \leq t'-1$. We proceed to prove this by induction. Note that by 
(\ref{eq42}) and fact that $\sigma_k \geq \sigma_{k-1} + 1$, it follows that (\ref{eq45}) holds 
for $r = t$. Suppose that (\ref{eq45}) holds for every $t \leq r \leq s$ for some $s$ which 
satisfies $t \leq s \leq t'-2$. Hence, $x_r^i \neq 0$ for every $i \in J$ and $t+1 \leq r \leq 
s$. It follows by Lemma \ref{lem16} that 
$$
\|{\bf x}_J^s - {\bf x}_J^*\| \leq \Delta \|{\bf x}_J^t - {\bf x}_J^*\| + \beta \max_{t \leq r \leq s} 
\|{\bf x}_{J^c}^r - {\bf x}_{J^c}^*\|_\infty + \mu \sum_{r=t}^{s-1} \|{\bf x}^r - 
{\bf x}^{r+1}\|^2. 
$$

\noindent
By (\ref{eq39}, (\ref{eq44}) and (\ref{eq45}), we obtain 
\begin{equation} \label{eq47}
\| {\bf x}_J^s - {\bf x}_J^*\| \leq \Delta \delta + \beta (\sigma_{k-1} \delta + \delta) + \mu 
\delta. 
\end{equation}

\noindent
Hence, for every $i \in J$, it follows from (\ref{eq40}), (\ref{eq42}), (\ref{eq44}), (\ref{eq47}), 
and the definition of $\sigma_k$ that 
\begin{eqnarray*}
|x_i^{s+1}| 
&\geq& |x_i^t| - \|{\bf x}_J^t - {\bf x}_J^{s+1}\|\\
&\geq& |x_i^t| - (\|{\bf x}_J^t - {\bf x}_J^*\| + \|{\bf x}_J^* - {\bf x}_J^s\| + \|{\bf x}_J^s - 
{\bf x}_J^{s+1}\|)\\
&>& \sigma_k \delta - (\delta + \Delta \delta + \beta \sigma_{k-1} \delta + \beta \delta 
+ \mu \delta + \delta)\\
&=& \sigma_{k-1} \delta + \delta. 
\end{eqnarray*}

\noindent
Thus, by induction, we conclude that (\ref{eq46}) holds for every $i \in J$ and $t \leq r \leq 
t'-1$. It follows by the arguments above that (\ref{eq47}) holds for every $t \leq s \leq t'-1$. 
Note that $\beta > 1$ and $\|{\bf y}\|_\infty \leq \|{\bf y}\|$ for any vector ${\bf y}$. 
Hence, by (\ref{eq45}) and the definition of $\sigma_k$, we obtain 
$$
\|{\bf x}^r - {\bf x}^*\|_\infty \leq (\Delta + \beta \sigma_{k-1} + \beta + \mu) \delta \leq 
\sigma_k \delta, 
$$

\noindent
for every $t \leq r \leq t'-1$. This proves part (b) of the required result. 

It follows from (\ref{eq40}) and (\ref{eq46}) with $r = t'-1$ that 
\begin{eqnarray*}
|x_i^{t'}| 
&\geq& |x_i^{t'-1}| - \|{\bf x}^{t'-1} - {\bf x}^{t'}\|\\
&>& \sigma_{k-1} \delta + \delta - \delta\\
&=& \sigma_{k-1} \delta. 
\end{eqnarray*}

\noindent
This proves part (a) of the required result. \hfill$\Box$ 

\bigskip

\noindent
{\it Proof of Lemma \ref{lem18}} We proceed by induction on $k$. If $k = 1$, then $J = I^*$. Hence, by 
(\ref{eq28}), $x_i^r = 0$ for every $i \notin J$ and $r \geq t$.  By Lemma \ref{lem17} (b), 
the claim holds for $k = 1$ (the proof of Lemma \ref{lem17} part (b) goes through 
verbatim even in $t' = \infty$).  Suppose now that the result holds for every $1 \leq k \leq 
h-1$, for some $h \geq 2$. Choose $J \subseteq I^*$ with $|J| \geq |I^*| - h + 1$ 
arbitrarily. Let $t > \max(r_0, r_1)$ be such that 
\begin{eqnarray}
& & |x_i^t| > \sigma_k \delta, \; \forall i \in J, \label{eq51}\\
& & |x_i^t| \leq \sigma_{k-1} \delta, \; \forall i \notin J. \label{eq52} 
\end{eqnarray}

\begin{enumerate}
\item {\bf Case 1}: $|x_i^r| \leq \sigma_{h-1} \delta$ for every $i \notin j$ and all $r 
\geq t$. 

\smallskip

\noindent
Since $|x_i^t| > \sigma_h \delta$ for every $i \in J$, it follows from Lemma \ref{lem17} 
part (b) that there exists an ${\bf x}^* \in \mathcal{X}^*$ such that 
$$
\|{\bf x}^r - {\bf x}^*\|_\infty \leq \sigma_h \delta, 
$$

\noindent
for every $r \geq t$. Hence the result holds for $k = h$, with $\bar{t} = t$. 
\item {\bf Case 2}: There exists an $r > t$ and $i \notin J$ such that $|x_i^r| > \sigma_{h-1} 
\delta$. 

\smallskip

\noindent
Let ${\bf t'}$ be the smallest $r > t$ such that $|x_i^r| > \sigma_{h-1} \delta$ for some 
$i \notin J$. By (\ref{eq52}), we obtain that $|x_i^r| \leq \sigma_{h-1} \delta$ for every 
$i \notin J$ and $t \leq r \leq t'-1$, and by (\ref{eq51}), $|x_i^t| > \sigma_h \delta$ for 
every $i \in J$. It follows by Lemma \ref{lem17} part (a) that $|x_i^{t'}| > \sigma_{h-1} 
\delta$ for every $i \in J$. 

Consider the $h+1$ intervals $T_0, T_2, \cdots, T_h$, where $T_0 = [0, \sigma_0 
\delta]$, $T_i = (\sigma_{i-1} \delta, \sigma_i \delta]$ for every $1 \leq i \leq h-1$, 
and $T_h = (\sigma_{h-1} \delta, \infty)$. Since $|x_i^{t'}| > \sigma_{h-1} \delta$ for 
some $i \notin J$, it follows that $T_h$ contains at least $|J| + 1$ entries (in absolute 
value) of the vector ${\bf x}^{t'}$. By (\ref{eq28}) and the fact that $\sigma_0 = 1$, we 
obtain that at least $n - |I^*|$ entries (in absolute value) of the vector ${\bf x}^{t'}$ are 
contained in $J_0$. Note that $|J| \geq |I^*| - h + 1$. Hence, this leaves at most $h - 2$ 
entries which are contained (in absolute value) in one of the $h-1$ intervals $T_1, \cdots, 
T_{h-1}$. By the Pigeon Hole principle, there exists a $q \in \{1,2, \cdots, h-1\}$ such that 
$|x_i^{t'}| \notin T_q$ for every $1 \leq i \leq n$. Let $h'$ denote the largest $q$ for which 
this occurs. Let $J' = \{j: |x_j^{t'}| > \sigma_{h'} \delta\}$. It follows from the observations 
above that $|J'| \geq |J| + h - h' \geq |J| + 1 - h'$. Note that by (\ref{eq28}), $J' \subseteq 
I^*$. Since $h' < h$, the induction hypothesis applied to $h'$, $t'$ and $J'$ yields the 
existence of an ${\bf x}^* \in \mathcal{X}^*$ and $\bar{t} \geq t'$ such that 
$$
\|{\bf x}^r - {\bf x}^*\| \leq \sigma_{h'} \delta, 
$$

\noindent
for every $r \geq \bar{t}$. Note that $\sigma_{h'} \leq \sigma_h$. Hence, the result holds 
for ${\bf x}^*$ and $k = h$. This completes the induction on $k$, and establishes the required 
result for every $1 \leq k \leq n$. 
\end{enumerate}

\noindent
\hfill$\Box$ 

\noindent
{\it Proof of Lemma \ref{lem19}} Fix any integer $\bar{r} \geq \max(r_0, r_1)$. 
\begin{enumerate}
\item {\bf Case 1}: $|x_i^r| \leq \sigma_n \delta$ for every $1 \leq i \leq n$ and $r \geq 
\bar{r}$. 

\smallskip

\noindent
In this case, let ${\bf x}^* \in \mathcal{X}^*$ be such that $\phi({\bf x}^{\bar{r}}) = 
\|{\bf x}^{\bar{r}} - {\bf x}^*\|$. It follows by (\ref{eq39}) that 
$$
|x_i^*| \leq |x_i^{\bar{r}}| + \|{\bf x}^{\bar{r}} - {\bf x}^*\| \leq \sigma_n \delta + 
\delta, 
$$

\noindent
for every $1 \leq i \leq n$. Note that $|x_i^r| \leq \sigma_n \delta$ for every $1 \leq i \leq 
n$ and $r \geq \bar{r}$, and by (\ref{eq28}), $x_i^r$ and $x_i^*$ are either both 
non-positive or non-negative. It follows that 
$$
\|{\bf x}^r - {\bf x}^*\|_\infty \leq \sigma_n \delta + \delta, 
$$

\noindent
for every $r \geq \bar{r}$. Hence, (\ref{eq53})holds with $\hat{r} = \bar{r}$. 
\item {\bf Case 2}: There exists $t \geq \bar{r}$ and $i \in \{1,2, \cdots, n\}$ such that 
$|x_i^t| > \sigma_n \delta$. 

\smallskip

\noindent
Let $\tilde{T}_j = (\sigma_{j-1} \delta, \sigma_j \delta]$ for $j = 1,2, \cdots, n$. Note that 
$|x_i^t|$ does not belong to any of $\tilde{T}_1,\tilde{T}_2, \cdots, \tilde{T}_n$. By the 
Pigeon Hole principle, there exists $q \in \{1,2, \cdots, n\}$ such that $T_q$ does not 
contain any entry (in absolute value) of the vector ${\bf x}^t$. Let $k$ be the largest $q$
for which this occurs. Let $\tilde{J}' = \{j: |x_j^t| > \sigma_k \delta\}$. Then $|J| \geq n - 
k + 1$, as $\tilde{T}_{k+1}, \tilde{T}_{k+2}, \cdots, \tilde{T}_n\}$ each contain at least 
one entry (in absolute value) of ${\bf x}^t$ and $|x_i^t| > \sigma_n \delta$. By the 
definition of $\tilde{J}'$, it follows that 
\begin{eqnarray*}
|x_j^t| > \sigma_k \delta, \; \forall j \in \tilde{J}',\\
|x_j^t| \leq \sigma_{k-1} \delta, \; \forall j \notin \tilde{J}'. 
\end{eqnarray*}

\noindent
It follows by (\ref{eq28}) and the fact that $\sigma_k \geq 1$ that $\tilde{J}' \subseteq 
I^*$. Hence, the assumptions of Lemma \ref{lem18} hold with $k$, $\tilde{J}'$ and $t$. 
It follows from Lemma \ref{lem18} that there exists an ${\bf x}^* \in \mathcal{X}^*$ and a 
$\bar{t} \geq t$ such that 
$$
\|{\bf x}^r - {\bf x}^*\|_\infty \leq \sigma_k \delta, 
$$

\noindent
for every $r \geq \bar{t}$. Since $\sigma_k \leq \sigma_n$, we conclude that (\ref{eq53}) 
holds with ${\bf x}^*$ and $\hat{r} = \bar{t}$. 
\end{enumerate}

\noindent
\hfill$\Box$ 

\bigskip

\noindent
{\it Proof of Lemma \ref{lemwlog5}} The result for $i \in S$ follows by repeating exactly the same arguments as in the 
proof of Lemma \ref{lem11} (b) for this case. The result for $i \in S^c$ follows by 
(\ref{eqwlog5}) and (\ref{eqwlog16}). \hfill$\Box$


\begin{thebibliography}{25}
\bibitem{Auslender:1981}
Auslender, A. (1978). Minimisation de fonctions localement lipschitziennes: applications à la 
programma- tion mi-convexe, mi-différentiable. In: {\it Mangasarian, O.L., Meyer, R.R., and 
Robinson, S.M. (eds.) Nonlinear Programming} {\bf 3}, pp. 429–460. Academic, New York. 

\bibitem{Bertsekas:Tsitsiklis:1989}
Bertsekas,D.P. and Tsitsiklis,J.N. (1989). {\it Parallel and Distributed Computation: Numerical 
Methods}, Prentice Hall, Englewood Cliffs, NJ, USA.

\bibitem{Bradley:Fayyad:1999}
Bradley, P.S., Fayyad, U.M. and Mangasarian, O.L. (1999). Mathematical programming for data mining: 
formulations and challenges, {\it INFORMS J. Comput.} {\bf 11}, 217–238. 

\bibitem{Burke:1985}
Burke, J.V. (1985). Descent methods for composite nondifferentiable optimization problems, 
{\it Math. Program.} {\bf 33}, 260–279. 

\bibitem{Chen:Donoho:1999}
Chen, S., Donoho, D. and Saunders, M. (1999). Atomic decomposition by basis pursuit, {\it SIAM J. 
Sci. Comput.} {\bf 20}, 33–61. 

\bibitem{Donoho:Johnstone:1994}
Donoho, D.L. and Johnstone, I.M. (1994). Ideal spatial adaptation by wavelet shrinkage, 
{\it Biometrika} {\bf 81}, 425– 455. 

\bibitem{Donoho:Johnstone:1995}
Donoho, D.L. and  Johnstone, I.M. (1995). Adapting to unknown smoothness via wavelet shrinkage, 
{\it J. Am. Stat. Assoc.} {\bf 90}, 1200–1224. 

\bibitem{Fletcher:1982}
Fletcher, R. (1982). A model algorithm for composite nondifferentiable optimization problems. 
{\it Math. Program. Study} {\bf 17}, 67–76.  

\bibitem{Friedman:Hastie:2007}
Friedman, J., Hastie, T. and Tibshirani, R. (2007). Sparse inverse covariance estimation with the 
graphical lasso, {\it Biostatistics} {\bf 9}, 432–441. 

\bibitem{Fu:1998}
Fu, W. (1998). Penalized regressions: the bridge vs the lasso, {\it Journal of Computational and 
Graphical Statistics} {\bf 3}, 397-416. 

\bibitem{Fukushima:Mine}
Fukushima, M. and Mine, H. (1981). A generalized proximal point algorithm for certain non-convex 
minimization problems, {\it Int. J. Syst. Sci.} {\bf 12}, 989–1000. 

\bibitem{FHT:2010}
Friedman, J., Hastie, T., and Tibshirani, R. (2010). Regularization paths for generalized linear 
models via coordinate descent, {\it Journal of Statistical Software} {\bf 33}. 

\bibitem{FHHT:2007}
Friedman, J., Hastie, T., Hofling, H. and Tibshirani, R. (2007). Pathwise coordinatewise optimization, 
{\it Annals of Applied Statistics} {\bf 1}, 302-332. 

\bibitem{Hoffman:1952}
Hoffman, A.J. (1952). On Approximate solutions of systems of linear inequalities, Journal 
of Research of the National Bureau of Standards {\bf 49}, 263-265.  


\bibitem{KOR:2014:CONCORD}
Khare, K., Oh, S. and Rajaratnam, B. (2014). A convex pseudo-likelihood framework for high 
dimensional partial correlation estimation with convergence guarantees, to appear in 
{\it Journal of the Royal Statistical Society B}. 


\bibitem{Kiwiel}
Kiwiel, K.C. (1986). A method for minimizing the sum of a convex function and a continuously 
differentiable function, {\it J. Optim. Theory Appl.} {\bf 48}, 437–449. 

\bibitem{KKB:2007}
Koh, K., Kim, S.-J. and Boyd, S. (2007). An interior-point method for large-scale 
$\ell_1$-regularized logistic regression, {\it J. Mach. Learn. Res.} {\bf 8}, 1519–1555. 

\bibitem{LLAN:2006}
Lee, S., Lee, H., Abeel, P. and Ng, A. (2006) Efficient l1-regularized logistic regression, In 
{\it Proceedings of the 21st National Conference on Artificial Intelligence, 2006}. 


\bibitem{Luo:Tseng:1:1989}
Luo, Z. and Tseng, P. (1989) On the convergence of a matrix splitting algorithm for the 
symmetric linear complementarity problem, Technical Report LIDS-P 1884, Laboratory 
for Information and Decision Systems, Massachusetts Institute of Technology. 

\bibitem{Luo:Tseng:1989}
Luo, Z. and Tseng, P. (1989) On the convergence of the coordinate descent method for 
convex differentiable minimization, Technical Report LIDS-P 1924, Laboratory for 
Information and Decision Systems, Massachusetts Institute of Technology. 

\bibitem{Luo:Tseng:1992}
Luo, Z. and Tseng, P. (1992) On the convergence of the coordinate descent method for 
convex differentiable minimization, {\it Journal of Optimization Theory and Applications} 
{\bf 72}, 7-35. 

\bibitem{Meinshausen:Buhlmann:2006}
Meinshausen, N. and Bu ̈hlmann, P. (2006). High-dimensional graphs and variable selection 
with the lasso, {\it The Annals of Statistics} {\bf 34}, 1436–1462. 

\bibitem{Park:Hastie:2007}
Park, M. and Hastie, T. (2007). $L_1$-regularization path algorithm for generalized linear 
models, {\it Journal of the Royal Statistical Society B} {\bf 69}, 659-677. 

\bibitem{Peng:2009}
Peng, J., Wang, P., Zhou, N., and Zhu, J. (2009). Partial correlation estimation by joint sparse 
regression models, {\it Journal of the American Statistical Association} {\bf 104}, 735– 746.

r

\bibitem{Richtarik:Takac:2011}
Richtarik, P. and Takac, M. (2011). Iteration complexity of randomized block-coordinate descent methods for minimizing a composite function, available on {\it Arxiv}. 

\bibitem{Yu:2008}
Rocha, G., Zhao, P., and Yu, B. (2008). A path following algorithm for sparse pseudo-likelihood 
inverse covariance estimation (SPLICE). Technical report, Statistics Department, UC Berkeley, 
Berkeley, CA. 

\bibitem{Saha:Tewari:2010}
Saha, A. and Tewari, A. (2010). On the Finite Time Convergence of Cyclic Coordinate Descent 
Methods, {\it CoRR} abs/1005.2146. 

\bibitem{Shevade:Keerthi:2003}
Shevade, S.K. and Keerthi, S.S. (2003). A simple and efficient algorithm for gene selection 
using sparse logistic regression, {\it Bioinformatics} {\bf 19}, 2246-2253. 

\bibitem{Sardy:Bruce:2001}
Sardy, S., Bruce, A. and Tseng, P. (2001). Robust wavelet denoising, {\it IEEE Trans. Signal Proc.} 
{\bf 49}, 1146–1152. 

\bibitem{Sardy:Tseng:2004}
Sardy, S. and Tseng, P. (2004). AMlet, RAMlet, and GAMlet: automatic nonlinear fitting of additive 
models, robust and generalized, with wavelets, {\it J. Comput. Graph. Stat.} {\bf 13}, 283–309. 

\bibitem{Tibshirani:1996}
Tibshirani, R. (1995). Regression selection and shrinkage via the lasso, {\it J. Royal Statist. Soc. 
B} {\bf 57}, pp. 267-288. 

\bibitem{Tibshirani:Saunders:2003}
Tibshirani, R., Saunders, M., Rosset, S. and Knight, K. (2005). Sparsity and smoothness via the 
fused lasso, {\it J. Royal Statist. Soc. B} {\bf 67}, 91-108. 

\bibitem{Tseng:2001}
Tseng, P. (2001). Convergence of block coordinate descent method for nondifferentiable 
minimization, {\it J. Optim. Theory Appl.} {\bf 109}, 473–492. 

\bibitem{Tseng:Yun:2009}
Tseng, P. and Yun, S. (2009). A coordinate gradient descent method for nonsmooth separable 
minimization, {\it Math. Program., Ser. B} {\bf 117}, 387–423. 

\bibitem{Yun:Toh:2011}
Yun, S. and Toh, K. (2011). A coordinate gradient descent method for $\ell_1$-regularized 
convex minimization, {\it Computational Optimization and Applications} {\bf 48}, 273-307. 
\end{thebibliography}
\end{document}